\documentclass[10pt]{amsart}
\usepackage[centertags]{amsmath}
\usepackage{amsfonts,mathrsfs}
\usepackage{newlfont}
\usepackage[all]{xy}
\usepackage{graphicx}
\usepackage[latin1]{inputenc}
\usepackage[english]{babel}
\usepackage[T1]{fontenc}
\usepackage{fancyhdr}
\usepackage{amscd}
\usepackage{amsthm}
\usepackage{amssymb}
\usepackage{verbatim}
\usepackage[a4paper,top=4.2cm,bottom=3.8cm,left=3.1cm,right=3.1cm]{geometry}
\newlength{\defbaselineskip}
\newcommand{\setlinespacing}[1]%
           {\setlength{\baselineskip}{#1 \defbaselineskip}}

\newtheorem{thm}{Theorem}[section]
\newtheorem{prop}[thm]{Proposition}

\newtheorem*{thmintro}{Theorem}

\newtheorem*{corintro}{Corollary}

\newtheorem{lem}[thm]{Lemma}
\newtheorem{cor}[thm]{Corollary}
\theoremstyle{definition}
\newtheorem{defn}[thm]{Definition}
\newtheorem{ex}[thm]{Example}
\theoremstyle{remark}
\newtheorem{rem}[thm]{Remark}

\date{}

\DeclareMathOperator{\Hom}{Hom} \DeclareMathOperator{\Ext}{Ext}

\newcommand \ap{\alpha_p}
\newcommand \an{\alpha_{p^n}}

\newcommand \bbf[1]{\textbf{\textit{#1}}}

\newcommand \g{\mathcal{G}}
\newcommand \glb{\g^{(\lambda)}}
\newcommand \gmu{\g^{(\mu)}}
\newcommand \gmup{\g^{(\mu^p)}}
\newcommand \glbp{\g^{(\lambda^{p})}}
\newcommand \glbn{\g^{(\lambda^{p^n})}}
\newcommand\gln{G_{\lambda,n}}
\newcommand \glx[1]{G_{{\lambda},{#1}}}
\newcommand \gmx[1]{G_{{\mu},{#1}}}
\newcommand \gm{\mathbb{G}_m}
\newcommand \Ga{\mathbb{G}_a}

\newcommand \h[3]{H^{#1}(#2,#3)}
\newcommand \hxg[1]{H^{1}(X,#1)}

\newcommand \lb{\lambda}

\newcommand \mup{\mu_p}
\newcommand \mun{\mu_{p^{n}}}

\newcommand \N{\mathbb{N}}

\newcommand \oo[1]{\mathcal{O}_{#1}}

\newcommand \Q{\mathbb{Q}}

\newcommand \Z{\mathbb{Z}}

\newcommand \ha{\mbox{$\hookrightarrow$}}
\newcommand \id{\operatorname{id}}

\newcommand \In{\subseteq}

\newcommand \lha{\hookrightarrow}
\newcommand \too{\longrightarrow}

\newcommand \mTo {\longmapsto}

\newcommand \on{\stackrel}
\renewcommand \phi{\mbox{$\varphi$}}

\newcommand \pt   {\otimes}

\renewcommand \rho{\mbox{$\varrho$}}

   \def\clE{{\cal E}}
 \def\clG{{\cal G}}

\newcommand{\RW[2]}{W_{#1}{\Omega^{#2}_{X}}}

\newcommand \fr{\operatorname{F}}

\def\scE{\mathscr{E}}

\DeclareMathOperator{\ext}{Ext}

\DeclareMathOperator{\Sp}{Spec}

\DeclareMathOperator{\Proj}{Proj}

\makeatletter
\renewcommand\l@subsection{\@tocline{2}{0pt}{5pc}{5pc}{}}
\makeatother

 


\begin{document}
\title[Models of $\Z/p^2\Z$]{Models of $\Z/p^2\Z$ over a d.v.r. of unequal characteristic}
\author{Dajano Tossici}
\begin{abstract}
Let $R$ be a discrete valuation ring  of unequal characteristic
which contains a primitive $p^2$-th root of unity. If $K$ is the
fraction field of $R$, it is well known that $(\Z/p^2 \Z)_K\simeq
\mu_{p^2,K}$. We prove that any finite and flat $R$-group scheme of
order $p^2$   isomorphic to $(\Z/p^2 \Z)_K$ on the generic fiber
(i.e. a model of $(\Z/p^2\Z)_K$), is the kernel in a short exact
sequence which generically coincides with the Kummer sequence. We
will explicitly describe and classify such models.
\end{abstract}
\address{Dipartimento di Matematica, Università di Roma Tre, Rome, Italy}
\email{dajano@mat.uniroma3.it} \maketitle
\begin{center}
\end{center}
\tableofcontents
\section*{Introduction}
 \textsc{Notation and Conventions.} 
If not otherwise specified we denote by $R$ a discrete valuation
ring (in the sequel d.v.r.) of unequal characteristic, i.e. a
d.v.r. with fraction field $K$ of characteristic zero and  with
residue field $k$ of characteristic $p>0$. Moreover we write
$S=\Sp(R)$. If, for $n\in \N$, there exists a distinguished
primitive $p^n$-th root of unity $\zeta_n$ in a d.v.r. $R$,
 we call $\lambda_{(n)}:=\zeta_n-1$. We remark that $v(\lb_{(n-1)})=p v(\lb_{(n)})$ and $ v(p)=p^{n-1}(p-1)v(\lb_{n})$. Moreover, for any $i\le n$, we suppose $\zeta_{i-1}=\zeta_i^p$.  And we will denote by
$\pi\in R$ one of its uniformizers. Moreover if $G$ is an affine
$R$-group scheme  we will denote by $R[G]$ the associated Hopf
algebra. All the schemes will be assumed no\oe therian.
%

\vspace{1cm}

Let $K$ be a field of characteristic $0$ which contains a
primitive $p^n$-th root of unity. We remark that this implies
$\mu_{p^n}\simeq\Z/p^n\Z$. We recall the following
 exact
sequence
$$
1\too \mu_{p^n}\too\gm\on{p^n}{\too} \gm\too 1,
$$
so-called the Kummer sequence. The Kummer theory says that any
$p^n$-cyclic Galois extension of $K$ can be deduced by the Kummer
sequence. We stress that the  Kummer sequence can be written also
as follows
$$
1\too \mu_{p^n}\too\gm^n\on{\theta_n}{\too} \gm^n\too 1
$$
where
$\theta_n((T_1,\dots,T_n))=(1-T_1^p,T_1-T_2^p,\dots,T_{n-1}-T_n^p)$.

Let $k$ be a field of characteristic $p>0$. The following exact
sequence
$$
0\too \Z/p^n\Z\too W_n(k)\on{\fr-1}\too W_n(k)\too 0,
$$
where $W_n(k)$ is the group scheme of Witt vectors of length $n$,
is called the Artin-Schreier-Witt sequence. The
Artin-Schreier-Witt theory implies that any $p^n$-cyclic Galois
covering of $k$ can be deduced by the Artin-Schreier-Witt
sequence.

 Let now $R$ be a d.v.r. of unequal characteristic which contains a $p^n$-th root of
unity. It has been proved, independently, by Oort-Sekiguchi-Suwa
(\cite{SOS}) and Waterhouse (\cite{Wat3})  the existence of an
exact sequence of group schemes over $R$ which unifies the above
two sequences for $n=1$. Later Green-Matignon (\cite{GM1}) and
Sekiguchi-Suwa(\cite{SS4}) have, independently, constructed
explicitly a unifying exact sequence for $n=2$.  This means that
it has been found an exact sequence
\begin{equation}\label{eq:KASWn}
0\too \Z/p^2\Z\too {\cal{W}}_{2}\too {\cal{W'}}_{2}\too 0
\end{equation}
 that  coincides with the Kummer sequence on the generic fiber and
with the Artin-Schreier-Witt sequence on the special fiber. The
case $n>2$  is treated in \cite{KASW2} and \cite{KASW1}. In this
paper we focus on finite and flat $R$-group schemes of order $p^2$
which are isomorphic to $(\Z/p^2\Z)_K$ on the generic fiber, i.e.
models of $(\Z/p^2\Z)_K$.  And we will prove that, for any such a
group scheme $G$, there exists an exact sequence
$$
0\too G\too \clE_1 \too \clE_2\too 0,
$$
with $\clE_1,\clE_2$ smooth $R$-group schemes, which coincides with
the Kummer sequence on the generic fiber. We will describe
explicitly all such isogenies and their kernels. Moreover we will
give a classification of models of $(\Z/p^2\Z)_K$.


 We now explain more precisely the classification we have
obtained. The first two sections are devoted to review some known
facts: in the first one we recall the definition of a class of
group schemes of order $p^n$, called $G_{\lb,n}$,  which are
isomorphic to $\mun$ on the generic fiber; in the second one we
recall some results about Neron blow-ups.


In the  third section we  recall the following well known result
about classification of models of $(\Z/p\Z)_K$.
\begin{thmintro}$\ref{teo:modelli di Z/pZ}$. Let suppose that $R$
contain a primitive $p$-{th} root of unity. If $G$ is a finite and
flat $R$-group scheme such that $G_K\simeq (\Z/p\Z)_K$ then $G\simeq
G_{\lb,1}$ for some $\lb\in R $.
\end{thmintro}

In \S \ref{sec:modelli di Z/p^2Z} we study the models of
$(\Z/p^2\Z)_K$. Let us suppose that $R$ contains a primitive
$p^2$-{th} root of unity. We remark that this hypothesis is only
used to conclude that $(\Z/p^2\Z)_K\simeq (\mu_{p^2})_K$. Without
this hypothesis what follows remains true substituting
$(\Z/p^2\Z)_K$ with $(\mu_{p^2})_K$.
 First of all  we show that any model of $(\Z/p^2\Z)_K$ is
an extension of $G_{\mu,1}$ by $G_{\lb,1}$ for some $\mu,\lb\in
R\setminus\{0\}$.
  So we  reduce ourselves to investigate on
$\Ext^1{(G_{\mu,1},G_{\lb,1})}$.

 Let $S_{\lb}:=\Sp(R/\lb R)$ and let us define the group
\begin{align*}
rad_{p,\lb}(<1+\mu S>):=\bigg\{&(F(S),j)\in
\Hom_{gr}({\gmx{1}}_{|S_\lb},{\gm}_{|S_\lb})\times \Z/p\Z \text{
such that }\\
& F(S)^p(1+\mu S)^{-j}=1\in
\Hom({\gmx{1}}_{|S_{\lb^p}},{\gm}_{|S_{{\lb^p}}})\bigg\}/<1+\mu
S>.
\end{align*}
There
 is an abuse  of notation since  $S$
 denotes both $\Sp(R)$ and the indeterminate of some polynomials. 
For any $(F,j)\in rad_{p,\lb}(<1+\mu S>)$ we will explicitly define
in \ref{subsec:explicit description...} an extension
$\clE^{(\mu,\lb;F,j)}$ of $G_{\mu,1}$ by $G_{\lb,1}$. As a group
scheme  $\clE^{(\mu,\lb;F,j)}$ is the kernel of an isogeny of smooth
group schemes of dimension $2$. This isogeny generically is
isomorphic to the morphism $\theta_n$ defined above.
 Using this notation, we will give a description of $\Ext^{1}(G_{\mu,1},\glx{1})$.
\begin{thmintro}$\ref{teo:ext1(glx,gmx)}$. Suppose that $\lb,\mu\in R$ with $v(\lb_{(1)})\ge v(\lb),v(\mu)$.
 There exists
an exact sequence
\begin{equation*}
\begin{array}{ll}
0\too rad_{p,\lb}(<1+\mu S>) \on{\beta}{\too}&
\Ext^{1}(G_{\mu,1},\glx{1})\too\\
&\too \ker \bigg(H^1(S,G_{\mu,1}^\vee)\too
H^1(S_\lb,G_{\mu,1}^\vee)\bigg) 
\end{array}
\end{equation*}
 where $\beta$ is defined by
$$
(F,j)\longmapsto \clE^{(\mu,\lb;F,j)}.
$$
In particular the set $\{\clE^{(\mu,\lb;F,j)}\}\In
\Ext^1(G_{\mu,1},\glx{1})$ is a group  isomorphic to
$rad_{p,\lb}(<1+\mu S>)$.
\end{thmintro}

The group $\Ext^{1}(G_{\mu,1},\glx{1})$   has been described by
Greither in \cite{gre} through  a short exact sequence, different by
that of the previous theorem. An advantage of our description is
that  we individuate a class of extensions which, if we forget the
structure of extension, "essentially" covers all the group schemes
of order $p^2$. Indeed from \ref{teo:ext1(glx,gmx)}  it follows that
any group scheme of order $p^2$, up to an extension of  $R$, is of
the form $\clE^{(\mu,\lb;F,j)}$ (see \ref{rem:essenzialmente tutti i
gruppi di ordine p2}).
Using the Sekiguchi-Suwa theory, which is briefly
explained in \S \ref{sec:Sekiguchi-Suwa theory}, we obtain the
following result.
\begin{corintro}$\ref{cor:clE se lb divide mu}$
Let us suppose $p>2$. Let $\mu, \lb \in R\setminus \{0\}$ be with
$v(\lb_{(1)})\ge v(\mu)\ge v(\lb)$. Then, the group $
\{\clE^{(\mu,\lb;F,j)} \} $  is isomorphic to the group
\begin{align*}
\Phi_{\mu,\lb}:=\bigg\{(a,j)\in {(R/\lb R)}\times\Z/p\Z& \text{
such that } a^p=0 \text{ and }pa-j\mu=\frac{p}{\mu^{p-1}}a^p\in
R/\lb^p R\bigg\},
\end{align*}
through the map
$$
(a,j)\longmapsto
\clE^{(\mu,\lb;\sum_{i=0}^{p-1}\frac{a^i}{i!}S^i,j)}.
$$
\end{corintro}
%
We  give also  a similar description of the group $
\{\clE^{(\mu,\lb;F,j)} \} $ with $v(\mu)<v(\lb)$: see
\ref{prop:rad_p}. Moreover we  remark that we can explicitly find
all the solutions $a$ of the equation $pa-j\mu\equiv
\frac{p}{\mu^{(p-1)}}a^p\mod \lb^{p}$  if $v(\mu)\ge v(\lb)$(see
\ref{cor:p2 surjective}).

In \S 5 we are interested  in the group schemes which are models of
$(\Z/p^2\Z)_K$. We prove the following theorem.
\begin{thmintro}$\ref{cor:modelli di Z/p^2 Z sono cosi'2}$
Let us suppose $p>2$. Let $G$ be a finite and flat $R$-group
scheme such that $G_K\simeq (\Z/p^2\Z)_K$. Then $G\simeq
\clE^{(\pi^{m},\pi^{n};\sum_{i=0}^{p-1}\frac{a^i}{i!}S^i,1)}$ for
some $v(\lb_{(1)})\ge m\ge n\ge 0$ and $(a,1)\in
\Phi_{\pi^m,\pi^n}$. Moreover $m,n$ and $a\in R/\pi^{n} R$ are
unique.
\end{thmintro}
The last section is devoted to determine, through the description of
\ref{cor:clE se lb divide mu}, the special fibers of the extensions
which, as group schemes, are models of $(\Z/p^2\Z)_K$.


The explicit description 
of the models of $(\Z/p^2\Z)_K$ presented in this paper will be
used in \cite{io3} to study the degeneration of $\Z/p^2\Z$-torsors
from characteristic $0$ to characteristic $p$.

\textbf{Acknowledgments} This paper constitutes part of my PhD
thesis. It is a pleasure to thank Professor Carlo Gasbarri for his
guidance and his constant encouragement. I  begun this work during
my visit to Professor Michel Matignon in 2005 at the Department of
Mathematics of the University of Bordeaux 1. I am indebted to him
for the useful conversations  and for his great interest in my
work.  I am deeply grateful  to Matthieu Romagny for  his very
careful reading of this paper and for his several comments,
suggestions, remarks and answers to my questions. I wish also to
thank Professor Fabrizio Andreatta for having suggested me the
sketch of proof of \ref{prop:torsori schemi normali}. Finally I
thank Filippo Viviani for the stimulating and useful
conversations.

\section{Some group schemes of order
$p^n$}\label{sec:glbn}

For any $\lambda\in R$ define the group scheme
$$
\glb=\Sp(R[T,\frac{1}{1+\lambda T}])
$$
The  $R$-group scheme  structure is given by
\begin{align*}
&T\too 1\pt T+T\pt 1 +\lb T\pt T\qquad&\text{comultiplication}\\
&T\too 0 \qquad& \text{counit}\\
&T\too -\frac{T}{1+\lb T}\qquad& \text{coinverse}
\end{align*}

We observe that if $\lb=0$ then $\glb\simeq \Ga$. It is possible
to prove that $\glb\simeq \gmu$ if and only if $v(\lambda)=v(\mu)$
and the isomorphism is given by $T\too \frac{\lambda}{\mu}T$.
Moreover it is easy to see that, if $\lb\in \pi R\setminus\{ 0\}$,
then $\glb_k\simeq \mathbb{G}_a$ and $\glb_K\simeq \gm$. It has
been proved by Waterhouse and Weisfeiler, in \cite[2.5]{Wat}, that
 any deformation, as a group scheme, of $\Ga$ to
 $\gm$ is isomorphic to $\glb$ for some $\lb\in \pi R\setminus\{0\}$.
If $\lb\in R\setminus\{0\}$ we can define the  morphism
\begin{align*}
\alpha^{\lb}:\glb\too \gm\\
\end{align*}
given, on the level of Hopf algebras,  by  $x\mTo 1+\lb x$: it is
an isomorphism on the generic fiber. If $v(\lb)=0$ then
$\alpha^{\lb}$ is an isomorphism.

We now define some finite and flat group schemes of order $p^n$.
Let $\lambda\in R$  satisfy the condition
 $$
(*)\qquad v(p)\ge p^{n-1}(p-1)v(\lambda).
$$ Then the map
\begin{align*}
\psi_{\lb,n}:&\glb\too \g^{(\lb^{p^n})}\\
             &T\too P_{\lb,n}(T):=\frac{(1+\lambda T)^{p^n}-1}{\lambda^{p^n}}
\end{align*}
is an isogeny of  degree $p^n$. Let
$$
\gln:=\Sp(R[T]/P_{\lb,n}(T))
$$
be  its kernel. It is a commutative finite flat group scheme over
$R$ of rank $p^n$. It is possible to prove that
\begin{align*}
{(G_{\lb,n})}_k\simeq  \mu_{p^n} &\qquad\text{ if $v(\lb)=0$}\\
{(G_{\lb,n})}_k\simeq \an &\qquad\text{ if $p^{n-1}(p-1)v(\lb)<v(p)$};\\
(G_{\lb,n})_k\simeq \alpha_{p^{n-1}}\times \Z/p\Z &\qquad \text{
if $p^{n-1}(p-1)v(\lb)=v(p)$}.
\end{align*}

  We observe that $\alpha^{\lambda}$ is compatible with $\psi_{\lb,n}$, i.e
 the following  diagram is commutative
\begin{equation}\label{eq:compatibilità phi_lb e alpha^lb}
\xymatrix@1{\glb\ar[d]_{\psi_{\lb,n}}  \ar[r]^{\alpha^{\lb}}&\ar[d]^{p^{n}}\gm\\
\glbn \ar[r]^(.6){\alpha^{\lb^{p^n}}}&\gm}
\end{equation}
Then it   induces a  map
$$
\alpha^{\lb,n}:\gln\too\mun
$$
which is an isomorphism on the generic fiber. And if $v(\lb)=0$
then $\alpha^{\lb,n}$ is an isomorphism.


We remark that
$$
 \Hom(G_{\lb,n},G_{\lb',n})=\left\{%
\begin{array}{ll}
    0, & \hbox{\text{ if } $ v(\lb)<v(\lb')$;} \\
    \Z/p^n\Z, & \hbox{\text{otherwise}.} \\
\end{array}%
\right.
$$
If $v(\lb)\ge v(\lb')$ the morphisms are given by
\begin{align*}
G_{\lb,n}&\too G_{\lb',n}\\
 T&\longmapsto \frac{(1+\lb T)^i-1}{\lb'}
\end{align*}
for $i=0,\dots,p^n-1$. It  follows easily that $G_{\lb,n}\simeq
G_{\lb',n}$ if and only if $v(\lb)=v(\lb')$.

 In the following any time we will speak
about $G_{\lb,n}$ it will be assumed that $\lb$ satisfies $(*)$. If
$R$ contains a primitive $p^n$-{th} root of unity $\zeta_n$ then,
since $$v(p)=p^{n-1}(p-1)v(\lb_{(n)}),$$  the condition $(*)$ is
equivalent to $v(\lb)\le v(\lb_{(n)})$.

\section{Néron blow-ups} 
 We recall here the definition of Néron blow-up. For details see
\cite[Ch. 3]{ner}, \cite{SS7} and \cite{Wat}. In this section  $R$
is a  not necessarily of unequal characteristic.
%
\begin{defn}
Let $X$ be a flat affine $R$-scheme of finite type and  $R[X]$ its
coordinate ring. Let $Y$ be a closed subscheme of $X_k$ defined by
some proper ideal $I(Y)$ of $R[X]$. Then $\pi\in I(Y)$. We define
the \textsl{Néron blow-up} (or \textsl{dilatation}) of $Y$ in $X$
by
$$
 X^{Y}:=\Sp(A[\pi^{-1}I(Y)]).
$$
\end{defn}
\vspace{0.5cm} Then  $X^Y$ is a flat affine $R$-scheme of finite
type and  the $R$-homomorphism $R[X]\In R[\pi^{-1}I(Y)]$ induces a
morphism
$$
X^Y\too X,
$$
which gives an isomorphism on the generic fiber.

 The Néron-blow up is explicitly given as follows: let
$I=(\pi,f_1,\dots,f_k)$ with $f_i\in R$. Then
$$
R[X^Y]=R[X][\pi^{-1}{f_1},\dots,\pi^{-1}f_k].
$$
So $X^Y$ is the open set of $x\in \Proj(\oplus_{i\ge 0}I^i)$ (the
classical blow-up of $X$ in $Y$), where $I_x$ is generated by
$\pi$. Clearly it is possible to give the definition for schemes
in general (see \cite[Ch. 3]{ner}).

In the following we are interested in the case where $X$ is an
affine
 flat group scheme $G$ and $Y$ a subgroupscheme $H$ of $G_k$. We  recall the following definitions.
\begin{defn} Let $\phi:G\too H$ be a morphism of flat $R$-group
schemes which is  an isomorphism restricted to the generic fibers.
Then it is called  a \textit{model map}.
\end{defn}
\begin{defn}
Let $H_K$ be a group scheme over $K$. Any  flat $R$-group scheme
$G$ such that $G_K\simeq H_K$ is called a \textsl{model} of $H_K$.
\end{defn}
It is possible to prove that $G^H$ is a group scheme and
$G^{H}\too G$ is a model map (\cite[1.1]{Wat}). We recall the
following results:
\begin{prop}\label{prop:prop. univ. blow-ups}
The canonical map $G^H\too G$ sends the special fiber into $H$.
Moreover $G^H$ has the following universal property: any model map
$G'\too G$ sending the special fiber into $H$ factors uniquely
through $G^H$.
\end{prop}
\begin{proof}
\cite[1.2]{Wat}.
\end{proof}
\begin{thm}\label{teo:decompositions blow-ups}
Any model map between affine group schemes is isomorphic to a
composite of Néron blow-ups.
\end{thm}
\begin{proof}
\cite[1.4]{Wat}.
\end{proof}
\begin{ex}\label{ex:blow-up glx}Let us consider the group scheme $G_{\mu,1}=\Sp(R[S]/(\frac{(1+\mu S)^p-1}{\mu^p}))$ with $v(p)>(p-1)v(\mu)$. The only possible subgroup of $(G_{\mu,1})_k$ which gives a nontrivial blow-up is $H=e$. Then $I(H)=(\pi,S)$ if $v(\mu)>0$ and $I(H)=(\pi,S-1)$ otherwise. It is easy to see
that, in both cases,
$$
\gmx{1}^e=G_{\mu\pi,1}.
$$
So if there exists a model map $G\too G_{\mu,1}$ then, using
\ref{teo:decompositions blow-ups}, $G\simeq G_{\lb,1}$ for some
$\lb\in R$.
\end{ex}

\section{Models of $(\Z/p\Z)_K$}

%
We now recall the classification of $(\Z/p\Z)_K$-models. Two
proofs of this result are for instance given in \cite[1.4.4,
3.2.2]{RoPhD}. The second one is essentially that we present here.
We remark that if $G$ is a model of $(\Z/m \Z)_K$ and $R$ contains
a primitive $m$-th root of unity then there are the following
model maps
$$
\Z/m\Z\too G\too \mu_m.
$$
Indeed the first one is the normalization map, while the second
one is the dual morphism of the normalization $\Z/m\Z\too
G^{\vee}$ (see also \cite[2.2.3]{Ra1} for a more general result).

\begin{prop}\label{teo:modelli di Z/pZ}
Let us suppose that $R$ contains a primitive $p$-{th} root of unity.
If $G$ is a finite and flat  $R$-group scheme such that $G_K\simeq
\Z/p\Z$ then $G\simeq G_{\lb,1}$ for some $\lb\in R\setminus\{0\}$.
\end{prop}
\begin{proof}
As  remarked  above we have an $R$-model map
$$
\phi:G\too \mu_p.
$$
By \ref{teo:decompositions blow-ups} it is a composition of Néron
blow-ups. Then, by \ref{ex:blow-up glx}, it follows that $G\simeq
G_{\lb,1}$ for some $\lb\in R\setminus\{0\}$.
\end{proof}
\section{Models of $(\Z/p^2\Z)_K$}\label{sec:modelli di Z/p^2Z}
In this section we study models of $(\Z/p^2\Z)_K$. Throughout the
section we  suppose that $R$ contains a primitive $p^2$-{th} root of
unity. First of all we prove that any such group is an extension of
$G_{\mu,1}$ by $G_{\lb,1}$ for some $\mu,\lb\in R\setminus\{0\}$.
\begin{lem}\label{lem:modelli di Z/p^2 Z sono estensioni}
Let $G$ be a finite and flat $R$-group scheme of order $p^2$ such
that $G_K$ is a constant group. Then $G$ is an extension of
$G_{\mu,1}$ by $G_{\lb,1}$ for some $\mu,\lb\in
R\setminus\{0\}$.\end{lem}
\begin{proof}If $G_K$ is a constant group then
$G_K$ is isomorphic to $(\Z/p^2\Z)_K$ or to $(\Z/p\Z)_K\times
(\Z/p\Z)_K$. We consider the factorization
$$
0\too (\Z/p\Z)_K\too G_K\too (\Z/p\Z)_K\too 0.
$$
We take the closure $G_1$ of $(\Z/p\Z)_K$ in $G$. Then $G_1$ is a
model of $(\Z/p\Z)_K$. So by \ref{teo:modelli di Z/pZ} it follows
that $G_1\simeq G_{\lb,1}$ for some $\lb\in R\setminus\{0\}$.
$G/G_{\lb,1}$ is a model of $(\Z/p\Z)_K$, too. So, again by
\ref{teo:modelli di Z/pZ}, we have $G/G_{\lb,1}\simeq G_{\mu,1}$
for some $\mu\in R\setminus\{0\}$. We are done.
\end{proof}
So we  study, first of all,  the group
$\Ext^1(G_{\mu,1},G_{\lb,1})$.
\subsection{Extensions of group schemes}

We here recall some generalities on extensions of group schemes.
For more details see \cite[III.6]{DG}.

Let $G$ and $H$ be  group schemes  on $S$. We moreover suppose
that $H$ is commutative and that $G$ acts on $H$. Let us denote
\begin{eqnarray*}
\Ext_S^0(G,H)=\{\phi\in \Hom_{Sc h_S}(G,H)|
\phi(gg')=\phi(g)+g(\phi(g')) \\
\text{ for any local sections}&\text{ $g,g'$ of $G$}\}.
\end{eqnarray*}
We are interested in  the case that $G$ acts trivially on $H$. In
this situation $$\Ext_S^0(G,H)=\Hom_{gr}(G,H).$$

Now $H\mapsto \Ext^{0}(G,H)$ is a left exact functor from the
category of fppf-sheaves of $G$-modules on $S$ to that of abelian
groups. Let $\Ext^{\bullet}_S(G,H)$ denote the left derived
functor of $H\mapsto \Ext^0_S(G,H)$. It is known that
$\Ext^1_S(G,H)$ is isomorphic to the group of equivalence classes
of extensions of $G$ by $H$ (see \cite[III 6.2]{DG}).

Recall that an extension of $G$ by $H$ is by definition an exact
sequence of fppf-sheaves of groups
$$
0\too H\on{i}{\too} E\on{j}{\too} G\too 0,
$$
such that $i(j(g)h)=gi(h)g^{-1}$ for any local sections $h$ of $H$
and $g$ of $E$.

Consider two extensions $(E):0\too H\on{i}{\too} E\on{j}{\too}
G\too 0$ and $(F):0\too H \on{i}{\too} F\on{j}{\too} G\too 0$.
They are equivalent if there exists a morphism of group schemes
$f:E\too F$ which makes the following diagram
$$
\xymatrix@1{(E):0\ar[r] &H\ar[r]^i\ar@{=}[d]& E\ar[r]^j\ar[d]^f
&G\ar@{=}[d]\ar[r] &0\\
(F):0\ar[r]& H\ar[r]^i& F\ar[r]^j& G\ar[r]& 0}
$$
commute. Clearly such  an $f$ is an isomorphism of group schemes.
If $G$ and $H$ are flat affine groups over $S$, then it is the
same for $E$.

We now recall the definitions of pushforward and pull-back of
extensions. Let $G$ and $H$ be  as above and $\phi:G'\too G$ a
morphism of group-schemes. Then $\phi$ induces a morphism
$$
\phi^*:\Ext^1_S(G,H)\too \Ext^1_S(G',H).
$$
It is explicitly  given as follows. Let
$$
(E):0\too H\on{i}{\too} E\on{j}{\too} G\too 0
$$
be an extension of $G$ by $H$. Then $\phi^*[E]$ is defined by the
diagram
$$
\xymatrix@1{\phi^*[E]:0\ar[r] &H\ar[r]^i\ar@{=}[d]&
E'\ar[r]^j\ar[d] &G'\ar[d]^\phi\ar[r]
&0\\
(E):0\ar[r]& H\ar[r]^i& E'\ar[r]^j& G\ar[r]& 0}
$$
where the right square is cartesian.

Now consider a group scheme $H'$ together with a $G$-action. If
$\psi:H\too H'$ is a morphism which preserves the $G$-action then
it induces a morphism
$$
\psi_*:\Ext^1_S(G,H)\too \Ext^1_S(G,H'),
$$
which we can explicitly describe as follows. Let
$$
(E):0\too H\on{i}{\too} E\on{j}{\too} G\too 0
$$
be an extension of $G$ by $H$. Then $\psi_*[E]$ is defined by the
diagram
$$
\xymatrix@1{(E):0\ar[r] &H\ar[r]^i\ar[d]^\psi& E\ar[r]^j\ar[d]
&G\ar@{=}[d]\ar[r]
&0\\
\psi_*[E]:0\ar[r]& H'\ar[r]^i& E'\ar[r]^j& G\ar[r]& 0}
$$
where the left square is cocartesian.

 Next we recall the
Hochschild cohomology. Let $G$ be a presheaf of groups on
$Sch_{|S}$ and $F$ a presheaf of $G$-modules on $Sch_{|S}$. We
define a complex $\{C^n(G,F),\delta^n\}$ as follows: $C^n(G,F)$
denotes the set of  morphisms of  schemes  from $G^n$ to $H$ and
the boundary map
$$
\delta^n:C^n(G,F)\too C^{n+1}(G,F)
$$
is defined by
\begin{eqnarray*}
(\delta^nf)(g_0,g_1,\dots,g_n)=g_0f(g_1,\dots,g_n)+\\
+\sum_{i=0}^{n-1}(-1)^{i+1}f(g_0,g_1,\dots,g_ig_{i+1},&\dots,g_n)+(-1)^{n+2}f(g_0,g_1,\dots,g_{n-1}).
\end{eqnarray*}

Put
\begin{gather*}
Z^n(G,F)=\ker (\delta^n:C^n(G,F)\too C^{n+1}(G,F)), \\
B^n(G,F)=Im(\delta^{n-1}:C^{n-1}(G,F)\too C^{n}(G,F)),
\end{gather*}
and
$$
H_0^n(G,F)=Z^n(F,G)/B^n(G,F).
$$

For our purposes we are interested in  the second group of
cohomology.  The following result  is indeed well  known.
\begin{prop}\label{prop:hochschild and extensions}
Let $G$ and $H$ be group schemes  over $S$. Given an action of $G$
on $H$ then $H_0^2(G,H)$ is isomorphic to the group of equivalence
classes of  extensions of $G$ by $H$ which have a scheme-theoretic
section.
\end{prop}
\begin{proof}
\cite{DG}.
\end{proof}
\subsection{Sekiguchi-Suwa Theory}\label{sec:Sekiguchi-Suwa theory}
Here is a very partial review  of  results of \cite{SS7},
\cite{SS6} and \cite{SS4}. Let $\mu,\lb\in \pi R\setminus\{0\}$.
For any $\lb'\in R\setminus\{0\}$  set $S_{\lb'}=\Sp(R/\lb' R)$.
 What we call Sekiguchi-Suwa theory is their description of
$\Hom_{gr}(\gmu_{|S_{\lb}},{\gm}_{|S_\lb})$ and
$\Ext^1(\gmu,\glb)$ through Witt vectors.

Let $Y=\Sp(R[T_1,\dots,T_m]/(F_1,\dots, F_n))$ be an affine
$R$-scheme of finite type. We recall that, for any $R$-scheme $X$ we
have that $\Hom_{Sch}(X,Y)$ is in bijective correspondence with the
set
$$\{(a_1,\dots,a_m)\in H^0(Y,\oo{Y})^m| F_1(a_1,\dots,
a_m)=0,\dots,F_n(a_1,\dots, a_m)=0 \}.$$  With an abuse of
notation we will identify these two sets. If $X$ and $Y$ are
$R$-group schemes
  we will also identify $\Hom_{gr}(X,Y)$ with a subset of $$\{(a_1,\dots,a_m)\in H^0(Y,\oo{Y})^m| F_1(a_1,\dots,
a_m)=0,\dots,F_n(a_1,\dots, a_m)=0 \}.$$

 We now fix presentations for the group schemes
$\gm$ and $\glb$  with $\lb \in \pi R$. Indeed we write
$\gm=\Sp(R[S,1/S])$ and $\glb=\Sp(R[S,1/1+\lb S])$. We remark again
that throughout the paper will be a conflict of notation since  $S$
will denote both $\Sp(R)$ and an indeterminate. But it should not
cause any problem. Before illustrating the Sekiguchi-Suwa theory we
see what happens when $\mu \in R^*$. In this case $\gmu\simeq \gm$,
and we have the following well known lemma.

\begin{lem} \label{lem:hom gm to gm}For any $\lb \in \pi R$ we have
$$
\Hom_{gr}({\gm}_{|S_{\lb}},{\gm}_{|S_\lb})=\{S^i\in R[S,1/S]|i\in
\Z\}.
$$
In particular if $v(\lb_{1})\ge v(\lb_{2})>0$, the restriction map
$$
\Hom_{gr}({\gm}_{|S_{{\lb}_1}},{\gm}_{|S_{\lb_{1}}})\too
\Hom_{gr}({\gm}_{|S_{{\lb}_2}},{\gm}_{|S_{\lb_{2}}})
$$
is an isomorphism.
\end{lem}
Moreover for the extensions group   we have
\begin{prop} For any $\lb \in \pi R\setminus\{0\}$, any
 $S$-action of $\gm$ on $\glb$ is trivial. Moreover
$$
\Ext^{1}(\gm,\glb)=0
 $$
 \end{prop}

 \begin{proof}
 See \cite[I 1.6, II 1.4]{SS6}.
 \end{proof}
 We also want to recall what happens to the extensions group when
 $\lb\in R^*$, i.e. $\glb\simeq \gm$. \begin{prop}
 For any $\mu \in R\setminus\{0\}$, any action of $\gmu$ on $\gm$
 is trivial. Moreover
 $$
\Ext^{1}(\gmu,\gm)=0.
 $$
 \end{prop}
 \begin{proof}
 See \cite[I 1.5, I 2.7]{SS6}
 \end{proof}

We now  consider the case $\mu,\lb\in \pi R\setminus\{0\}$. Any action of $\gmu$ on $\glb$ is trivial (\cite[I 1.6]{SS6}). For any flat $R$-scheme $X$ let us consider 
 the exact sequence on the fppf site $X_{fl}$
\begin{equation}\label{eq:succ esatta glb in gm}
0\too \glb\on{\alpha^{\lb}}{\too}{\gm}\too
i_*{{\mathbb{G}}_{m,X_{\lb}}}\too 0,
\end{equation}
 where $i$ denotes the closed immersion $X_\lb=X\pt_R(R/\lb R)\ha X$
(see \cite[1.2]{SS2}). We observe that by definitions we have that
$$
\Hom_{gr}(\gmu_{|S_{\lb}},{\gm}_{|S_\lb})=\{F(S)\in (R/\lb
R[S,\frac{1}{1+\mu S}])^*|F(S)F(T)=F(S+T+\mu ST) \}
$$
 If
we apply the functor $\Hom(\cdot,\gm)$ to the sequence
\eqref{eq:succ esatta glb in gm} we obtain, in particular,  a
 map 
\begin{equation*}
\begin{aligned}
\Hom_{gr}(\gmu_{|S_\lb},{\gm}_{|S_\lb})
\on{\alpha}{\too}\Ext^1(\gmu&,\glb).
\end{aligned}
\end{equation*}
 given by
$$
F\longmapsto \clE^{(\mu,\lb;F)},
$$
where
$$
\clE^{(\mu,\lb;F)}
$$
is a smooth affine  commutative group defined as follows: let
$\tilde{F}(S)\in R[S]$ be a lifting of $F(S)$, then
$$
\clE^{(\mu,\lb;F)}= \Sp(R[S_1,S_2,\frac{1}{1+\mu
S_1},\frac{1}{\tilde{F}(S_1)+\lb S_2}])
$$
\begin{enumerate}
    \item law of multiplication
    \begin{align*}
    S_1\longmapsto &S_1\pt 1+1\pt S_1+\mu S_1\pt S_1\\
    S_2\longmapsto &S_2\pt \tilde{F}(S_1)+\tilde{F}(S_1)\pt S_2 +\lb S_2\pt S_2+\\
                   & \qquad   \quad    \frac{\tilde{F}(S_1)\pt \tilde{F}(S_1)-\tilde{F}(S_1\pt 1+1\pt S_1+\mu S_1\pt S_1)}{\lb}
    \end{align*}
    \item unit
    \begin{align*}
    &S_1\longmapsto 0\\
    &S_2\longmapsto \frac{1-\tilde{F}(0)}{\lb}
    \end{align*}
    \item inverse
    \begin{align*}
    &S_1\longmapsto -\frac{S_1}{1+\mu S_1}\\
& S_2 \longmapsto\frac{\frac{1}{\tilde{F}(S_1)+\lb
S_2}-\tilde{F}(-\frac{S_1}{1+\mu S_1})}{\lb}
\end{align*}
\end{enumerate}
We moreover  define the following  homomorphisms of group schemes
$$
\glb=\Sp(R[S,(1+\lb S)^{-1}])\too \clE^{(\mu,\lb;F)}
$$
by
\begin{align*}
&S_1\longmapsto 0\\
&S_2\longmapsto S +\frac{1-\tilde{F}(0)}{\lb}
\end{align*}
and
$$
\clE^{(\mu,\lb;F)}\too \gmu=\Sp(R[S,\frac{1}{1+\mu S}])
$$
 by
\begin{align*}
&S\too S_1.
\end{align*}
It is easy to see that
\begin{equation}\label{eq:estensioni lisce}
0\too \glb\too \clE^{(\mu,\lb;F)}\too \gmu\too 0
\end{equation}
is exact. A different choice of the lifting $\tilde{F}(S)$ gives
an isomorphic extension. We recall the following theorem.
\begin{thm}\label{teo:ss ext1}For any $\lb,\mu\in \pi R\setminus\{0\}$, the map $$\alpha:\Hom_{gr}(\gmu_{|S_\lb},{\gm}_{|S_\lb}){\too}
\Ext^1(\gmu,\glb)$$  is a surjective morphism of groups. And
$\ker(\alpha)$ is generated by the class of $1+\mu S$. In
particular any extension of $\gmu$ by $\glb$ is commutative.
\end{thm}
\begin{proof}
\cite[\S 3]{SS9}.
\end{proof}
 \noindent We now define some spaces which had been used by Sekiguchi and Suwa to describe $\Hom_{gr}(\gmu_{|S_\lb},{\gm}_{|S_\lb})$ and, by the above result,  $\Ext^1(G_{\mu,1},G_{\lb,1})$.  See \cite{SS4} for details.
\begin{defn}
For any ring $A$, let $W_n(A)$ be  the ring of Witt vectors of
length $n$ and $W(A)$  the ring of infinite Witt vectors. We
define

\begin{align*}
\widehat{W}_n(A)=\bigg\{(a_0,\dots,a_n)\in W_n(A)|& a_i\text{ is
nilpotent for any } i  \text{ and } \\
&a_i=0 \text{ for all but a finite number of } i \bigg\}
\end{align*}
and
\begin{align*}
\widehat{W}(A)=\bigg\{(a_0,\dots,a_n,\dots)\in W(A)|& a_i\text{ is
nilpotent for any } i  \text{ and } \\
&a_i=0 \text{ for all but a finite number of } i \bigg\}.
\end{align*}
\end{defn}

We recall the definition of the so-called Witt-polynomial: for any
$r\ge 0$ it is
$$
\Phi_r(T_0,\dots,T_r)=T_0^{p^r}+pT_{1}^{p^{r-1}}+\dots+p^rT_r.
$$
Then  the following maps are defined:
\begin{itemize}
    \item[-]Verschiebung
    \begin{align*}
    V:W_{n}(A)&\too W_{n+1}(A)\\
(a_0,\dots,a_n)&\longmapsto (0,a_0,\dots,a_n)
    \end{align*}
    \item[-]Generalization  of Frobenius
 \begin{align*}
    \fr:W_{n+1}(A)&\too W_{n}(A)\\
(a_0,\dots,a_n)&\longmapsto
(F_0(\mathbf{T}),F_1(\mathbf{T}),\dots,F_n(\mathbf{T}))
    \end{align*}
where the polynomials $F_r(\mathbf{T})=F_r(T_0,\dots,T_r)\in
\Q[T_0,\dots,T_{r+1}]$ are defined inductively by
$$
\Phi_r(F_0(\mathbf{T}),F_1(\mathbf{T}),\dots,F_r(\mathbf{T}))=\Phi_{r+1}(T_0,\dots,T_{r+1}).
$$
\end{itemize}
If $p=0\in A$ then $\fr$ is the usual Frobenius. The subring
$\widehat{W}(A)$  is stable respect to these maps.

For any morphism $G: \widehat{W}(A)\too \widehat{W}(A)$ we will
set $\widehat{W}(A)^G:=\ker G$. And for any $a\in A$ we denote
 the element $(a,0,0,\dots,0,\dots)\in W(A)$ by $[a]$.

We recall the following standard result about Witt vectors.
\begin{lem}\label{lem:isobaricita della somma
}Let $S_r[\mathbf{T},\mathbf{U}]\in \Z[\mathbf{T},\mathbf{U}]$
such that, if $\bbf{a},\bbf{b}\in W(A)$, then
$$
\bbf{a}+\bbf{b}=(S_0[\bbf{a},\bbf{b}],\dots,S_r[\bbf{a},\bbf{b}],\dots)
$$
If $T_i$ and $U_i$ have weight $p^i$ then
$S_r[\mathbf{T},\mathbf{U}]$ is isobaric of weight $p^r$.

\end{lem}

The following lemma  will be useful later.
\begin{lem}\label{lem:somma termine per termine} Let $\lb \in R$. If
$\bbf{a}=(a_0,a_1,\dots),\bbf{b}=(b_0,b_1,\dots)\in
\widehat{W}(R/\lb R)^{\fr}$ then
$$
\bbf{a}+\bbf{b}=(a_0+b_0,a_1+b_1,\dots, a_i+b_i,\dots)
$$
\end{lem}
\begin{proof}
We suppose that $\bbf{a}+\bbf{b}=(c_0,c_1,\dots,c_i,\dots)$. By
the previous lemma we have that $c_r(\bbf{a},\bbf{b})$ is isobaric
of weight $p^r$. It is a standard result that
$$c_r(\bbf{a},\bbf{b})=a_r+b_r+
c'_r((a_0,a_1,\dots,a_{r-1}),(b_0,b_1,\dots,b_{r-1})).$$ for some
polynomial $c'_r(S_0,\dots, S_{r-1},T_0,\dots,T_{r-1})$. Clearly
$c'_r(\bbf{a},\bbf{b})$ is isobaric of weight $p^r$, too. Hence
$\deg (c'_r)\ge p$.

Let $\tilde{a_i},\tilde{b_i}\in R$ be liftings of $a_i$ and $b_i$,
respectively. For any $r\ge 1$, up to changing $\bbf{a}$ with
$\bbf{b}$, we can suppose that
$v(\tilde{a_k})=\min\{v(\tilde{a_i}),v(\tilde{b_i})|
i=0,\dots,r-1\}$, for some $0\le k \le r-1$. Since $\deg c'_r\ge p$
then $v( c'_r(\tilde{\bbf{a}},\tilde{\bbf{{b}}}))\ge
pv(\tilde{a}_k)$. But $v(\tilde{a}_k^p)\ge v(\lb)$ since
$\fr(\bbf{a})=0$. Hence $c'_r(\bbf{a},\bbf{b})=0\in R/\lb R$. So
$$
\bbf{a}+\bbf{b}=(a_0+b_0,a_1+b_1,\dots, a_i+b_i,\dots)
$$
\end{proof}
  We now recall the definition of the
Artin-Hasse exponential series
$$
E_p(T):=\exp\bigg(\sum_{r\ge
0}\frac{T^{p^r}}{p^r}\bigg)=\prod_{r=0}^{\infty}\exp\bigg(\frac{T^{p^r}}{p^r}\bigg)\in
\Z_{(p)}[[T]].
$$
Sekiguchi and Suwa introduced  a deformation of the Artin-Hasse
exponential map in \cite{SS4}. By the well known formula
$\lim_{\mu\to 0 }(1+\mu x)^{\frac{\alpha}{\mu}}=\exp(\alpha x)$,
it can be seen that $(1+\mu x)^{\frac{\alpha}{\mu}}$ is a
deformation of $\exp(\alpha x)$. From this point of view they
defined the formal power series $E_p(U,\Lambda;T)\in
\Q[U,\Lambda][[T]]$ by
$$
E_p(U,\Lambda;T):=(1+\Lambda
T)^{\frac{U}{\Lambda}}\prod_{r=1}^{\infty}(1+\Lambda^{p^r}T^{p^r})^{\frac{1}{p^r}((\frac{U}{\Lambda})^{p^r}-(\frac{U}{\Lambda})^{p^{r-1}})}
$$
They proved that $E_p(U,\Lambda;T)$ has in fact its coefficients
in $\Z_{(p)}[U,\Lambda]$. It is possible to show (\cite[2.4]{SS4})
that
$$E_p(U,\Lambda;T)=\left\{%
\begin{array}{ll}
    \prod_{(i,p)=1}E_p(U\Lambda^{i-1}T^i)^\frac{(-1)^{i-1}}{i}, & \hbox{if $p>2$;} \\
    \prod_{(i,2)=1}E_p(U\Lambda^{i-1}T^i)^\frac{1}{i}\bigg[\prod_{(i,2)=1}E_p(U\Lambda^{2i-1}T^{2i})^\frac{1}{i}\bigg]^{-1}, & \hbox{if $p=2$.} \\
\end{array}%
\right.
$$
Let $A$ be a $\Z_{(p)}$-algebra and $a,\mu\in A$. We define
$E_p(a,\mu;T)$ as $E_p(U,\Lambda;T)$ evaluated at $U=a$ and
$\Lambda=\mu$.
\begin{ex}\label{es:esempi di Ep(a,mu,T)} It is easy to see that  $E_p(a,0;T)=E_p(a
T)$ and $E_p(\mu,\mu;T)=1+\mu T$. Moreover if $a^p=\mu^{p-1}a\in
A$ then $(\frac{a}{\mu})^{p^r}-(\frac{a}{\mu})^{p^{r-1}}=0$ for
$r\ge 1$. Hence
$$
E_p(a,\mu;T)=(1+\mu
T)^{\frac{a}{\mu}}=1+\sum_{i=1}^{p-1}\frac{\prod_{k=0}^{i-1}(a-k\mu)}{i!}T^i.
$$
In particular if $\mu=0$ and $a^p=0\in A$ then
$$
E_p(a,0;T)=\sum_{i=0}^{p-1}\frac{a^i}{i!}T^i.
$$
\end{ex}

\vspace{.5cm}

If \textbf{\textit{a}}$=(a_0,a_1,a_2,\dots)\in W(A)$ we define the
formal power series
\begin{equation}
E_p(\textbf{\textit{a}},\mu;T)=\prod_{k=0}^\infty
E_p(a_k,\mu^{p^k};T^{p^k}).
\end{equation}

 The following  result gives an explicit description of
 $\Hom_{gr}(\gmu_{|A},{\gm}_{|A})$.

\begin{thm}\label{teo:ss hom} Let $A$ be a
$\Z_{(p)}$-algebra and $\mu\in A$ a nilpotent element. 
The homomorphism
\begin{align*}
\xi^0_A:\widehat{W}(A)^{\fr-[\mu^{p-1}]}&\too \Hom_{gr}(\gmu_{|A},{\gm}_{|A})\\
\textbf{a}&\longmapsto E_p(\textbf{a},\mu;S)
\end{align*}
is bijective.
\end{thm}
\begin{proof}
\cite[2.19.1]{SS4}.
\end{proof}

And \ref{teo:ss ext1} and \ref{teo:ss hom} give the following:
\begin{cor} For any $\lb,\mu \in \pi R\setminus\{0\}$ the map\begin{align*}
\alpha\circ \xi_{R/\lb R}^0 :\widehat{W}(R/\lb R)^{\fr-[\mu^{p-1}]}/<1+\mu T>&\too \ext^1(\gmu,{\gm})\\
\textbf{a}&\longmapsto \clE^{(\lb,\mu;E_p(\textbf{a},\lb;S))}
\end{align*}
is an isomorphism.
\end{cor}

We now describe some natural maps through these identifications.
Consider the isogeny
$$
\psi_{\mu,1}:\gmu\too \gmup.
$$
Let us now suppose that $p>2$. Then we have that,  if $p^2\equiv 0
\mod \lb$,
$$
\psi_{\mu,1}^*:\Hom_{gr}({\gmup}_{|S_{\lb}},{\gm}_{|S_{\lb}})\too
\Hom_{gr}(\gmu_{|S_{\lb}},{\gm}_{|S_{\lb}})
$$
is given by
\begin{gather}\label{eq:phi*:W(R)too W(R)}
\textit{\textbf{a}}\longmapsto[\frac{p}{\mu^{p-1}}]\textit{\textbf{a}}+V(\textit{\textbf{a}})
\end{gather}
(see \cite[1.4.1 and 3.8]{SS4}).

For $p=2$ the situation is slightly different. Let us define
 a variant of the Verschiebung as follows. Define polynomials
 $$
\tilde{V}_r(\mathbf{T})=\tilde{V}_r(T_0,\dots,T_r)\in
\Q[T_0,\dots,T_{r}]
$$
 inductively  by $\tilde{V}_0=0$ and
$$
\Phi_r(\tilde{V}_0(\mathbf{T}),\dots,\tilde{V}_r(\mathbf{T}
))=p^{p^r}\Phi_{r-1}(T_0,\dots,T_{r-1})
$$
for $r\ge 1$. Then we have that (with possibly $2^2\not\equiv
0\mod \lb$)
$$
\psi_{\mu,1}^*:\Hom_{gr}(\clG^{(\mu^2)}_{|S_{\lb}},{\gm}_{|S_{\lb}})\too
\Hom_{gr}(\gmu_{|S_{\lb}},{\gm}_{|S_{\lb}})
$$
is given by
\begin{gather*}
\textit{\textbf{a}}\longmapsto[\frac{2}{\mu}]\textit{\textbf{a}}+V(\textit{\textbf{a}})+\tilde{V}(\textit{\textbf{a}})
\end{gather*}
(see \cite[ 3.8]{SS4}).

For simplicity, to avoid to use this description of
$\psi_{\mu,1}^*$,  we will consider sometimes  only the case $p>
2$. 

 Consider the morphism
\begin{equation}\label{eq:def p}
\begin{aligned}
 p:\Hom_{gr}(\gmu_{|S_{\lb}},{\gm}_{|S_\lb})&\too\Hom_{gr}(\gmu_{|S_{\lb}{^p}},{\gm}_{|S_{\lb^p}})\\
     F(S)&\longmapsto F(S)^p
     \end{aligned}
\end{equation}

This morphism  is such that
$$
{\psi_{\lb,1}}_*\circ \alpha=\alpha\circ p.
$$
Let \textbf{\textit{a}}$\in (\widehat{W}(R/\lb
R))^{\fr-[\mu^{p-1}]}$. Take any lifting
$\tilde{\textbf{\textit{a}}}\in W(R)$. Using the identifications
of \ref{teo:ss hom} the morphism $p$ above
is given by
\begin{equation}\label{eq:p a}
\textbf{\textit{a}}\longmapsto p\tilde{\textbf{\textit{a}}}
\end{equation}
(see \cite[4.6]{SS4}).  We will sometimes simply write $p\bbf{a}$.

\subsection{Two exact sequences}\label{sec:two exact sequences}
The main tools which we will use to 
calculate the extensions of $\glx{1}$ by $G_{\mu,1}$ are two exact
sequences. We recall them in this subsection. See \eqref{eq:ker
phi} and \eqref{eq:ker alpha} below. First of all we prove that
any action of $G_{\mu,1}$ on $G_{\lb,1}$ is trivial.
\begin{lem}\label{lem:morfismo gen nullo=nullo}Let $\phi:G\too H$  be an $S$-morphism of affine
$S$-groups. Assume that $G$ is flat over $S$. Then $\phi=0$ if and
only if the generic fiber $\phi_K=0$.
\end{lem}
\begin{proof}
\cite[1.1]{SS6}.
\end{proof}

\begin{lem}
Every action of $\gmx{1}$  on $\glx{1}$ is trivial.
\end{lem}
\begin{proof}
Giving an action of $\glx{1}$ on $\gmx{1}$ is the same as giving a
morphism $\gmx{1}\too Aut_R(\glx{1})$. If we consider the generic
fiber we have a morphism
$$
\mu_{p,K}\too Aut_{K}(\mu_{p,K}).
$$
The last one is the étale group scheme $(\Z/p \Z)^*_K$. It is a
group scheme of order $p-1$. So any morphism $\mu_{p,K}\too
Aut_K(\mu_{p,K})$ is trivial. Applying \ref{lem:morfismo gen
nullo=nullo} we have the thesis.
\end{proof}

In the following, all the actions will be supposed trivial.
Applying now the functor $\Ext$ to the following exact sequence of
group schemes
$$
(\Lambda): \qquad 0\too
\glx{1}\on{i}{\too}\glb\on{\psi_{\lb,1}}{\too} \g^{(\lb^{p})}\too
0,
$$
we obtain
\begin{equation}
\begin{aligned}\label{eq:ker phi}
0\too\Hom_{gr}(G_{\mu,1},\glbp)\on{\delta'}{\too}
\Ext^{1}(G_{\mu,1},\glx{1})\on{i_*}\too& \\
\too \Ext^{1}(G_{\mu,1},\glb)
&\on{\psi_{{\lb,1}_*}}{\too}\Ext^{1}(\gmx{1},\glbp).
\end{aligned}
\end{equation}
We remark that $\delta'$ is injective since
$$
{\psi_{\lb,1}}_*:\Hom_{gr}(G_{\mu,1},\glb)\too
\Hom_{gr}(G_{\mu,1},\glbp)
$$
is the zero morphism.
 Indeed since $G$ is flat over $R$, then by \ref{lem:morfismo gen
nullo=nullo},
$$
\Hom_{gr}(\gmx{1},\glb)\lha \Hom_{gr}
(\mu_{p,K},{{\gm}_{K}})\simeq\Z/p\Z.
$$
And it is easy to verify that
\begin{equation}\label{eq:Hom(gmu1,gl)}
\Hom_{gr}(\gmx{1},\glb)=\left \{%
\begin{array}{ll}
   \Z/p\Z, & \hbox{if $\lb\mid\mu$} \\
    0, & \hbox{if $\lb \nmid \mu$} \\
\end{array}%
\right.
\end{equation}
Let us write $\gmx{1}=\Sp(R[S]/(\frac{(1+\mu S)^p-1}{\mu^p}))$. If
$\lb \mid \mu$ the group is formed by the morphisms given by
$\sigma_i:S\longmapsto \frac{(1+\mu S)^{i}-1}{\lb}$ with $i\in
\Z/p\Z$.  The map $(\psi_{\lb,1})_*:\Hom_{gr}(G_{\mu,1},\glb)\too
\Hom_{gr}(G_{\mu,1},\glbp)$ is moreover nothing else but the
multiplication by $p$. So  it is clearly zero.

The map
\begin{equation}\label{eq:def delta'}
\delta':\Hom_{gr}(G_{\mu,1},\glbp)\too \Ext^{1}(G_{\mu,1},\glx{1})
\end{equation}
is defined by
$$
\sigma_i\longmapsto (\sigma_i)^*(\Lambda),
$$
where $(\sigma_i)^*(\Lambda)$ is explicitly
$$
\Sp(R[S_1,S_2]/(\frac{(1+\mu S_1)^p-1}{\mu^p},\frac{(1+\lb
S_2)^p-(1+\mu S_1)^i}{\lb^p}),
$$
with the maps
\begin{align*}
\glx{1}&\too \sigma_i^*(\Lambda)\\
   S_1&\longmapsto 0\\
   S_2&\longmapsto S
\end{align*}
and
\begin{align*}
 \sigma_i^*(\Lambda)&\too\gmx{1}\\
   S&\longmapsto S_1
\end{align*}
The structure of group scheme on $ \sigma_i^*(\Lambda)$ is the
unique one which makes the map
\begin{align*}
\sigma_i^*(\Lambda)&\too
\mu_{p^2}=\Sp(R[Z_1,Z_2]/(Z_1^p-1,Z_2^p-Z_1^i))\\
         Z_1  &\longmapsto  1+\mu S_1\\
         Z_2  &\longmapsto       1+\lb S_2
\end{align*}
a morphism of group schemes.

As remarked in \cite[4.4]{SS6}, there is the following long exact
sequence
\begin{equation}\label{eq:ker alpha}
\begin{array}{ll}
0\too &\Hom_{gr}(\gmx{1},\g^{(\lb')})\too
\Hom_{gr}(\gmx{1},\gm)\on{r_{\lb'}}{\too}
\Hom_{gr}({\gmx{1}}_{|S_{\lb'}},{\gm}_{|S_{\lb'}})\on{\delta}{\too}\\
&\too \Ext^{1}(\gmx{1},\g^{(\lb')})\on{\alpha^{\lb'}_{*}}{\too}
\Ext^{1}(\gmx{1},\gm)\too
\Ext^{1}({\gmx{1}}_{|S_{\lb'}},{\gm}_{|S_{\lb'}}).
\end{array}
\end{equation}

 We  so have
\begin{equation}\label{eq:ker alpha con r lambda}
\ker{\alpha^{\lb'}_{*}}\simeq Im \delta\simeq
\Hom_{gr}({\gmx{1}}_{|S_{\lb'}},{\gm}_{|S_\lb'})/r_{\lb'}(\Hom_{gr}(\gmx{1},\gm)).
\end{equation}

We remark that by \eqref{eq:Hom(gmu1,gl)}, setting $\lb'=1$, it
follows that
$$
\Hom_{gr}(\gmx{1},\gm)\simeq\Z/p\Z
$$
and the group is formed by the morphisms $S\longmapsto (1+\mu
S)^i$. While, by \eqref{eq:Hom(gmu1,gl)},
$\Hom_{gr}(\gmx{1},\g^{(\lb')})\simeq \Z/p\Z$ if  $\lb' |\mu$ and
it is $0$ otherwise. Hence, by \eqref{eq:ker alpha}, if $\lb'|\mu
$ then $r_{\lb'}$ is the zero morphism, otherwise $r_{\lb'}$ is an
isomorphism. 
Hence, by \eqref{eq:ker alpha con r lambda},
\begin{equation}\label{eq:ker(alpha)}
 \ker{\alpha^{\lb'}_{*}}\simeq
\Hom_{gr}({\gmx{1}}_{|S_{\lb'}},{\gm}_{|S_{\lb'}})/<1+\mu S>.
\end{equation}

 In the following we give a more explicit description of the main
ingredients of the exact sequences \eqref{eq:ker phi} and \eqref{eq:ker alpha}.

\subsection{Explicit description of
$\Hom_{gr}({\gmx{1}}_{|S_{\lb}},{\gm}_{|S_{\lb}})$}\label{subsec:explicit
description...}

First we consider the simplest cases. If $\lb \in \pi R$,
\begin{equation}\label{eq:Hom(mup,gm)su Slb}
\Hom_{gr}({\mup}_{|S_\lb},{\gm}_{|S_\lb})=\{S^i\in (R/\lb
R)[S,1/S]|i\in \Z/p\Z\}.
\end{equation}
While if $\lb\in R^*$ we have $S_\lb=\emptyset$ and
$\Hom_{gr}({\gmx{1}}_{|S_\lb},{\gm}_{|S_\lb})=\{   1 \}$.

Now we study $\Hom_{gr}({\gmx{1}}_{|S_{\lb}},{\gm}_{|S_{\lb}})$
for $\mu,\lb \in \pi R\setminus\{0\}$.
\begin{prop}\label{lem:suriettività mappa tra hom}
Let    $\lb,\mu \in \pi R\setminus\{0\}$. 
The map
$$
i^*:\Hom_{gr}(\gmu_{|S_\lb},{\gm}_{|S_\lb})\too
\Hom_{gr}({G_{\mu,1}}_{|S_\lb},{\gm}_{|S_\lb})
$$
induced by
$$
i:G_{\mu,1}\hookrightarrow \gmu
$$
is surjective. If $p>2$, $\xi_{R/\lb}^0$, defined in \ref{teo:ss
hom}, induces an isomorphism
\begin{align*}
\widehat{W}(R/\lb R)^{\fr-[\mu^{p-1}]}/\bigg\{
[\frac{p}{\mu^{p-1}}]\textit{\textbf{b}}+ V(\textit{\textbf{b}})
|\bbf{b} \in {\widehat{W}(R/\lb
R)}^{\fr-[{\mu^{p(p-1)}}]}\bigg\}\too
\Hom_{gr}({G_{\mu,1}}_{|S_\lb},{\gm}_{|S_\lb})
\end{align*}
\end{prop}
\begin{proof}
We have  by definitions that
\begin{align*}
\Hom_{gr}({\gmx{1}}_{|S_{\lb}},{\gm}_{|S_{\lb}})=\{F(S)\in
\bigg(R/\lb R[S]/(&\frac{(1+\mu S)^p-1}{\mu^p})\bigg)^*|\\
&F(S)F(T)=F(S+T+\mu ST)\}
\end{align*}
and
\begin{align*}
\Hom_{gr}(\gmu_{|S_{\lb}},{\gm}_{|S_{\lb}})=\{F(S)\in &\bigg(R/\lb
R[S,\frac{1}{1+\mu S}]\bigg)^*|\\ &F(S)F(T)=F(S+T+\mu ST)\}.
\end{align*}
 Since $(G_{\mu,1})_k$ is isomorphic to $\alpha_p$ or
$\Z/p\Z$ then the group $ \Hom_{gr}((\gmx{1})_k,{\gm}_{,k})$ is
trivial. So $F(S)\equiv 1\mod \pi$. Moreover any $F(S)\in (R/\lb
R)[S]/(\frac{(1+\mu S)^p-1}{\mu^p})$ such that $F(S)\equiv 1\mod
\pi$ is invertible. The same is true in
$\Hom_{gr}({\gmu}_{|S_\lb},{\gm}_{|S_\lb})$ since $\gmu_k\simeq
\Ga$. We now say that $F$ satisfies condition $(\sharp)$ if
\begin{itemize}
  \item[] $F(S)\equiv 1 \mod \pi  $;
  \item[] $F(S)F(T)=F(S+T+\mu  ST)$.
\end{itemize}
Then
$$
\Hom_{gr}({\gmx{1}}_{|S_{\lb}},{\gm}_{|S_{\lb}})=\{F(S)\in R/\lb
R[S]/\big(\frac{(1+\mu S)^p-1}{\mu^p}\big)| F(S) \text { satifies
} (\sharp)\}
$$
and
$$
\Hom_{gr}(\gmu_{|S_{\lb}},{\gm}_{|S_{\lb}})=\{F(S)\in R/\lb
R[S,\frac{1}{1+\mu S}]| F(S) \text { satifies } (\sharp)\}.
$$
 Any $F\in
R/\lb R[S]/(\frac{(1+\mu S)^p-1}{\mu^p})$ can be represented by a
polynomial of degree $p-1$. And if it satisfies  $(\sharp)$, it
also satisfies $(\sharp)$  in $R/\lb R[S,\frac{1}{1+\mu S}]$. 

So
$$
i^*:\Hom_{gr}(\gmu_{|S_\lb},{\gm}_{|S_\lb})\too
\Hom_{gr}({G_{\mu,1}}_{|S_\lb},{\gm}_{|S_\lb})
$$
is surjective.

 Now, by the exact sequence
$$
(\Lambda')  \qquad 0\too G_{\mu,1}\on{i}{\too}
\gmu\on{\psi_{\mu,1}}{\too} \gmup\too 0
$$
over $S_\lb$, we have the long exact sequence of cohomology
\begin{equation}\nonumber\label{eq:succ. esatta per hom  verso gm}
\begin{aligned}
0\too
\Hom_{gr}(\gmup_{|S_\lb},{\gm}_{|S_\lb})\on{\psi_{\mu,1}^*}{\too}
\Hom_{gr}(\gmu_{|S_\lb},{\gm}_{|S_\lb})\on{i^*}{\too}\\
\too
\Hom_{gr}({G_{\mu,1}}_{|S_\lb},{\gm}_{|S_\lb})\on{\delta''}{\too}&
\Ext^{1}(\gmup_{|S_\lb},{\gm}_{|S_\lb})\too \dots
\end{aligned}
\end{equation}
By \ref{teo:ss hom} we have that
$$
\Hom_{gr}(\gmu_{|S_\lb},{\gm}_{|S_\lb})\simeq \widehat{W}(R/\lb
R)^{F-[\mu^{p-1}]}
$$
and, by \eqref{eq:phi*:W(R)too W(R)},
$$
\psi_{\mu,1}^*(\Hom_{gr}(\gmup_{|S_{\lb}},{\gm}_{|S_{\lb}}))\simeq
\bigg\{ [\frac{p}{\mu^{p-1}}]\textit{\textbf{b}}+
V(\textit{\textbf{b}}) |\bbf{b} \in {\widehat{W}(R/\lb
R)}^{\fr-[{\mu^{p(p-1)}}]}\bigg\}.
$$
Therefore the proposition is proved.
\end{proof}

%

We now give a more explicit description  of
$\Hom_{gr}({\gmx{1}}_{|S_{\lb}},{\gm}_{|S_\lb})$.
\begin{prop}\label{pro:Hom(gmu,gm)}
If $\lb,\mu\in R$ with $v(p)\ge(p-1)v(\mu)>0$ and $v(p)\ge
v(\lb)$, then
\begin{align*}
\Hom_{gr}({\gmx{1}}_{|S_{\lb}},{\gm}_{|S_\lb})=\{E_p(a,\mu;
S)=1+\sum_{i=1}^{p-1}&\frac{\prod_{k=0}^{i-1}(a-k\mu)}{i!}S^i\\
 \text{s.t.
}a&\in R/\lb R \text{ and } a^p=\mu^{p-1} a\in R/\lb R\}
\end{align*}

\end{prop}

\begin{rem}

In \cite[3.5]{SS7}, an inductive formula for the coefficients of
the polynomials $F(T)\in \Hom({\gmu}_{|S_\lb},{\gm}_{|S_\lb})$ is
given. If we consider only polynomials of degree less or equal to
$p-1$, it coincides with \eqref{eq:formula induttiva}.  But for
the reader's convenience, we prefer to give here a direct proof of
this formula.
\end{rem}
\begin{rem}\label{rem:v(mu)ge v(lb)}If $v(\mu)\ge v(\lb)$ then
$$
\Hom_{gr}({\gmx{1}}_{|S_\lb},{\gm}_{|S_\lb})=\{\sum_{i=0}^{p-1}\frac{a^i}{i!}T^i|
a^p=0\}=\{E_p(aT)|a^p=0\}
$$

\end{rem}
\begin{proof}
%
As seen in \ref{lem:suriettività mappa tra hom}
\begin{align*}
\Hom_{gr}({\gmx{1}}_{|S_{\lb}},{\gm}_{|S_{\lb}})=\{F(S)=\sum_{i=0}^{p-1}{a_i}S^i\in
&R/\lb R[S]/(\frac{(1+\mu S)^p-1}{\mu^p})\\
\text{s.t. } F(S)\equiv 1 \mod \pi  &\text{ and }
F(S)F(T)=F(S+T+\mu ST)\}.
\end{align*}

Now
\begin{equation}\label{eq:F in termini di a}
\begin{aligned}
F(S+T+\mu ST)&=\sum_{i=0}^{p-1}a_i(S+T+\mu ST)^i\\
            &=\sum_{i=0}^{p-1}\sum_{j=0}^{i}\sum_{k=0}^{j}\binom{i}{j}\binom{j}{k}\mu^{i-j}a_iS^{k+i-j}T^{i-k}\\
            &=\sum_{r=0}^{p-1}\sum_{l=0}^{p-1}\sum_{\max\{r,l\}\le i\le
            r+l}\binom{i}{2i-(r+l)}\binom{2i-(r+l)}{i-l}\mu^{r+l-i}a_i
            S^rT^l
\end{aligned}
\end{equation}
and
\begin{equation}\label{eq:F(S)F(T)}
F(S)F(T)=\sum_{r=0}^{p-1}\sum_{l=0}^{p-1}a_ra_lS^rT^l.
\end{equation}
So we have the equality if and only if
\begin{equation}\label{eq:formula induttiva}
\begin{aligned}
a_ra_l&=\sum_{\max\{r,l\}\le i\le
r+l}\binom{i}{2i-(r+l)}\binom{2i-(r+l)}{i-l}\mu^{r+l-i}a_i \\
      &=\sum_{\max\{r,l\}\le i\le
r+l}\frac{i!}{(r+l-i)!(i-l)!(i-r)!}\mu^{r+l-i}a_i
\end{aligned}
\end{equation} for any $0\le r,l\le p-1$.
Clearly $a_0= 1$.

 We now have the following lemma:
\begin{lem}\label{lem:equivalenza statements}
For any $\mu,\lb\in \pi R\setminus \{0\}$, the following
statements are equivalent
\begin{itemize}
  \item [i)]$a_{r}= \frac{\prod_{k=0}^{r-1}(a_1-k\mu)}{r!}$ for any $1 \le r\le p-1$ and   $
 {\prod_{k=0}^{p-1}(a_1-k\mu)}= 0
    $;
    \item[ii)]$a_{r-1}a_1= (r-1)\mu a_{r-1}+ra_{r}
$ for any $1\le r\le {p-1}$;
  \item [iii)]$a_ra_l=\sum_{\max\{r,l\}\le i\le
r+l}\frac{i!}{(r+l-i)!(i-l)!(i-r)!}\mu^{r+l-i}a_i$ for any $1\le
l,r\le {p-1}$.
\end{itemize}
\end{lem}

\begin{proof}In the following we use the convention that $a_i=0$ if $i>p-1$.

$i)\Leftrightarrow ii)$.  It is clear that
$$
a_{r-1}a_1= (r-1)\mu a_{r-1}+ra_{r}
$$
 is equivalent to  $a_r= a_{r-1}\frac{a_1-\mu
(r-1)}{r}$, if $r<p$, and $a_{p-1}({a_1-\mu (p-1)})= 0$. An easy
induction shows that this is equivalent to
$$
a_{r}= \frac{\prod_{k=0}^{r-1}(a_1-k\mu)}{r!}
    $$
   if $r<p-1$ and
    $$
 {\prod_{k=0}^{p-1}(a_1-k\mu)}= 0.
    $$

$ii) \Leftarrow iii)$. It is obvious.

$ii)\Rightarrow iii)$. We will prove it by induction on $l$. By
hypothesis $a_{r-1}a_1= (r-1)\mu a_{r-1}+ra_{r}$ for any $r$. We
now suppose  that $iii)$ is true for $k\le l-1$ for any $r$. Then
we will prove it is also true for $l$ for any $r$. We can clearly
suppose $l\le r$, otherwise, up to a change of $l$ with $r$, we
can conclude by induction. We have
$$
\begin{array}{lll}
a_r a_l&\on{ii)}{=}\displaystyle{a_r(a_{l-1}\frac{a_1-\mu(l-1)}{l}) }\\
      &=\displaystyle{(a_ra_{l-1})\frac{a_1-\mu(l-1)}{l}}&\\
      &\on{induct.}{=}\displaystyle{\bigg(\sum_{r\le i\le r+l-1}\frac{i!}{(r+l-i-1)!(i-l-1)!(i-r)!}\mu^{r+l-i-1}a_i\bigg)\frac{a_1-\mu(l-1)}{l}}\\
     &\on{induct.}{=}\displaystyle{\sum_{r\le i\le
     r+l-1}\frac{i!}{(r+l-i-1)!(i-l-1)!(i-r)!}\mu^{r+l-i-1}}\\
     &\qquad\displaystyle{(\frac{\mu(i-l+1)a_i+(i+1)a_{i+1}}{l})}\\
     &= \displaystyle{\frac{r!}{l!(r+1-l)!}\mu^l(r+1-l)a_{r}+}\\
     &\qquad\displaystyle{+\sum_{r+1\le i\le
     r+l-1}\bigg(\frac{i!}{(r+l-i)!(i-l)!(i-r-1)!l}\mu^{r+l-i}
     + }\\&\qquad\displaystyle{+
     \frac{i!(i-l+1)}{(r+l-i)!(i-l+1)!(i-r)!l}\mu^{r+l-i}\bigg)a_i}+\frac{(r+l-1)!(r+l)}{r!l!}a_{r+l}\\
     &\displaystyle{= \sum_{r\le i\le
r+l}\frac{i!}{(r+l-i)!(i-l)!(i-r)!}\mu^{r+l-i}a_i}.
\end{array}
$$

\end{proof}

We come back to the proof of  the proposition. In $R/\lb R$  the
condition $${\prod_{k=0}^{p-1}(a_1-k\mu)}= 0
    $$ is equivalent to $a_1^p= \mu^{p-1} a_1$. Indeed we have the following equality in $\Z/p\Z[S]$
$$
\prod_{k=0}^{p-1} (S-k)=S^p-S,
$$
since these polynomials have the same zeros. Since $p= 0\in R/\lb
R$, then
$$
\prod_{k=0}^{p-1}(a_1-k\mu)= a_1^p-\mu^{p-1} a_1.
$$

By the lemma and \ref{es:esempi di Ep(a,mu,T)} the thesis follows.
\end{proof}
We now essentially rewrite \ref{lem:equivalenza statements} in a
more expressive 
form.
\begin{cor}\label{cor:equivalenze per F}
Let $\lb,\mu\in \pi R\setminus\{0\}$ and let
$F(S)=\sum_{i=0}^{p-1}a_iS^i\in R/\lb R[S]$ be a polynomial of
degree less than or equal to $p-1$. Then the following statements
are equivalent
\begin{itemize}
\item[(i)] $F(S)F(T)-a_0^2=F(S+T+\mu  ST)-a_0$\item[(ii)]
$F(S)a_1= F'(S)(1+\mu S) $
 where $F'$ is the formal derivative of $F$.
\end{itemize}
\end{cor}
\begin{rem}\label{rem:soluzione equazione diff}Let us suppose $v(\mu)\ge v(\lb)$. This corollary, together with \ref{pro:Hom(gmu,gm)}, says
that the solution of the differential equation in $R/\lb
R[S]/(\frac{(1+\mu S)^p-1}{\mu^p}) $
$$
\left\{%
\begin{array}{ll}
   F'(S)= a F(S) ,  \\
    F(0)= 1 \\
\end{array}%
\right.
$$
has as unique solution $F(S)=E_p(a
S)=\sum_{i=0}^{p-1}\frac{a^i}{i!}S^i$ and  $a^p= 0$.
\end{rem}
\begin{proof}By \eqref{eq:formula induttiva},  we have that  \ref{lem:equivalenza
statements}(iii) is equivalent to
\begin{equation}\label{eq:equazione}
F(S)F(T)-a_0^2=F(S+T+\mu  ST)-a_0.
\end{equation}
If we put $l=1$  in \eqref{eq:formula induttiva}, we obtain the
coefficient of $T$ in both members of \eqref{eq:equazione}. This
means that \ref{lem:equivalenza statements}(ii) is equivalent to
$$
F(S)a_1= F'(S)(1+\mu S).
$$
 Then their equivalence comes from \ref{lem:equivalenza
statements}.

\end{proof}
\noindent When $v(\mu)\ge v(\lb)$, putting together
\ref{lem:suriettività mappa tra hom} and \ref{pro:Hom(gmu,gm)}, we
have a simpler description of
$\Hom_{gr}({G_{\mu,1}}_{|S_\lb},{\gm}_{|S_\lb})$.
\begin{cor}\label{cor: hom gmu gm se lb divide mu}Let $p>2$.
Let  $\lb,\mu\in R$ with $v(p)\ge(p-1)v(\mu)>0$ and $v(\mu)\ge
v(\lb)>0$. Then we have the following isomorphism of groups
$$
(\xi^0_{R/\lb R})_p:(R/\lb R)^{\fr}
\too\Hom_{gr}({\gmx{1}}_{|S_{\lb}},{\gm}_{|S_\lb})
$$
given by
$$
a\longmapsto E_p(a S).
$$
Moreover the restriction map
$$
i^*:\Hom_{gr}(\gmu_{|S_{\lb}},{\gm}_{|S_\lb})\simeq\widehat{W}(R/\lb
R)^{\fr}\too
\Hom_{gr}({\gmx{1}}_{|S_{\lb}},{\gm}_{|S_\lb})\simeq(R/\lb
R)^{\fr}
$$
is given, in terms of Witt vectors,  by
$$
\bbf{a}=(a_0,a_1,\dots,0,0,0,\dots)\longmapsto
\sum_{i=0}^{\infty}(-1)^i(\frac{p}{\mu^{p-1}})^ia_i
$$
\end{cor}
\begin{proof}
We first remark that the restriction of the Teichm\"{u}ller map
$$
T:(R/\lb R)^{\fr}\too \widehat{W}(R/\lb R)^{\fr},
$$
given by
$$
a\longmapsto [a],
$$
is a morphism of groups. This follows from \ref{lem:somma termine
per termine}.
Moreover, if we consider the isomorphism $$\xi^0_{R/\lb
R}:\widehat{W}(R/\lb R)^{\fr}\too
\Hom_{gr}(\gmu_{|S_{\lb}},{\gm}_{|S_\lb})$$ and
$$i^*:\Hom_{gr}(\gmu_{|S_{\lb}},{\gm}_{|S_\lb})\too
\Hom_{gr}({\gmx{1}}_{|S_{\lb}},{\gm}_{|S_\lb}),$$ we have
$$
i^*\circ \xi^0_{R/\lb R}\circ T=(\xi^0_{R/\lb R})_p.
$$
So $(\xi^0_{R/\lb R})_p$ is a morphism of groups. It is surjective
by \ref{pro:Hom(gmu,gm)} and, by \ref{lem:suriettività mappa tra
hom}, its kernel is
$$T((R/\lb R)^{\fr})\cap \bigg\{
[\frac{p}{\mu^{p-1}}]\textit{\textbf{b}}+ V(\textit{\textbf{b}})
|\bbf{b} \in {\widehat{W}(R/\lb R)}^{\fr}\bigg\}.$$ Let us now
suppose that there exists $\bbf{b}=(b_0,b_1,\dots)\in
{\widehat{W}(R/\lb R)}^{\fr}$ and $ a\in (R/\lb R)^{\fr}$ such
that $[\frac{p}{\mu^{p-1}}]\textit{\textbf{b}}+
V(\textit{\textbf{b}})=[a]$.
 It follows by the definition of Witt vector ring that
\begin{equation}\label{eq:p mu^(p-1) b=...b_0 I}
[\frac{p}{\mu^{p-1}}]\textit{\textbf{b}}=(\frac{p}{\mu^{p-1}}b_0,
\dots, (\frac{p}{\mu^{p-1}})^{p^j} b_j,\dots),\end{equation} and
\begin{equation}\label{eq:a-V(b)}
[a]-V(\textit{\textbf{b}})=(a_0,-b_0,-b_1,\dots).
\end{equation}


Since $\bbf{b} \in {\widehat{W}(R/\lb R)} $, there exists $r\ge 0$
such that $b_j=0$ for any $j\ge r$. Moreover, comparing
\eqref{eq:p mu^(p-1) b=...b_0 I} and \eqref{eq:a-V(b)} it follows
\begin{align*}
&(\frac{p}{\mu^{p-1}})^{p^j}b_{j+1}=-b_j \quad\text{ for } j\ge 0\\
&(\frac{p}{\mu^{p-1}})^pb_0= a.
\end{align*}
Hence $b_j=a=0$ for any $j\ge 0$. It follows that $(\xi^0_{R/\lb
R})_p$ is injective.

We now prove the second part of the statement. First of all we
remark that for any $\bbf{a}=(a_0, \dots, a_j, \dots)\in
\widehat{W}(R/\lb R)^{\fr}$ we have
$$
\bbf{a}=\sum_{j=0}^{\infty}V^j([a_j]).
$$
It is clear that for any $a\in R/\lb R$ we have $i^*([a])=a$.
While, by \ref{lem:suriettività mappa tra hom}, it follows that
$i^*{V(\bbf{b})}=-i^*([\frac{p}{\mu^{p-1}}]\bbf{b})$ for any
$\bbf{b}\in \widehat{W}(R/\lb R)^{\fr}$. Hence
$i^*{V^j(\bbf{b})}=(-1)^ji^*([(\frac{p}{\mu^{p-1}})^j]\bbf{b})$
for any $j\ge 1$. From these facts it follows that
\begin{align*}
i^*(\bbf{a})&=i^*(\sum_{j=0}^{\infty}V^j([a_j]))\\
            &=\sum_{j=0}^{\infty}(i^*(V^j([a_j])))\\
            &=\sum_{j=0}^{\infty}(-1)^j(\frac{p}{\mu^{p-1}})^ja_j.
\end{align*}
\end{proof}
    \subsection{Explicit description of $\delta$}\label{par:delta}
The map
$$
\delta:\Hom_{gr}({\gmx{1}}_{|S_{\lb}},{\gm}_{|S_\lb})\too
\Ext^{1}(\gmx{1},\g^{(\lb)})
$$
can also be explicitly described. 
  We  have the following
commutative diagram
$$
\xymatrix{\Hom_{gr}(\gmu_{|S_\lb},{\gm}_{|S_\lb})\ar[r]\ar[d]^{\alpha}&
\Hom_{gr}({G_{\mu,1}}_{|S_\lb},{\gm}_{|S_\lb})\ar[d]^{\delta}\ar[r]&0\\
\Ext^{1}(\gmu,\g^{(\lb)})\ar[r]^{i^*}&\Ext^{1}(\gmx{1},\g^{(\lb)})}
$$
where the first horizontal map is surjective by
\ref{lem:suriettività mappa tra hom}. So, given  $$F(S)\in
\Hom_{gr}({\gmx{1}}_{|S_{\lb}},{\gm}_{|S_\lb}),$$ we can choose a
representant in  $\Hom_{gr}(\gmu_{|S_{\lb}},{\gm}_{|S_\lb})$ which
we denote again by $F(S)$ for simplicity. Then $\delta$ is defined
by
$$
F(S)\longmapsto
\widetilde{\clE}^{(\mu,\lb;F)}:=i^*(\clE^{(\mu,\lb;F)})=i^*(\alpha(F(S))).
$$

 If $\tilde{F}(S)\in R[S]$ is any lifting then it is
defined, as  a scheme, by
$$
\widetilde{\clE}^{(\mu,\lb;F)}=\Sp\big(R[S_1,S_2,(\tilde{F}(S_1)+\lb
S_2)^{-1}]/\frac{(1+\mu S_1)^p-1}{\mu^p}\big).
$$
%
This extension does not depend on the choice of the lifting since
the same is true for ${\clE}^{(\mu,\lb;F)}$.

So, by \eqref{eq:ker alpha}, we see that $\ker(\alpha^\lb_*)\In
\Ext^{1}(\gmx{1},\glb)$ is nothing else but the group $
\{\widetilde{\clE}^{(\mu,\lb;F)} \}.
$ We recall that, by \eqref{eq:ker alpha con r lambda}, for any
$\lb'\in R$,
\begin{equation*}
\ker{\alpha^{\lb'}_{*}}\simeq
\Hom_{gr}({\gmx{1}}_{|S_{\lb'}},{\gm}_{|S_{\lb'}})/r_{\lb'}(\Hom_{gr}(\gmx{1},\gm)). 
\end{equation*}  
We have therefore proved the following proposition.
\begin{prop}\label{prop:ker alpha=clE}Let  $\lb,\mu\in R\setminus\{0\}$ with $v(p)\ge(p-1)v(\mu)$.
Then $\delta$ induces an isomorphism 
\begin{align*}
 \Hom_{gr}({\gmx{1}}_{|S_{\lb}},{\gm}_{|S_\lb})/r_{\lb'}(\Hom_{gr}(\gmx{1},\gm))&\too \{\widetilde{\clE}^{(\mu,\lb;F)} \}\\
F(S)&\longmapsto \widetilde{\clE}^{(\mu,\lb;F)}
\end{align*}
\end{prop}

\begin{rem}\label{rem:i unici.}
As remarked in \S \ref{sec:two exact sequences}, if $\lb'\mid \mu$
then
$$r_{\lb'}(\Hom_{gr}(\gmx{1},\gm))=0,$$ otherwise
$$r_{\lb'}(\Hom_{gr}(\gmx{1},\gm))\simeq<1+\mu S>\simeq \Z/p\Z.$$
\end{rem}
\begin{ex}

\begin{itemize}

\item[a)] Let us suppose $v(\mu)=0$ and $v(\lb)>0$. Since, by
\eqref{eq:Hom(mup,gm)su Slb} and the previous remark,
$$\Hom_{gr}({\mu_p}_{|S_\lb},{\gm}_{|S_{\lb}})\simeq
r_{\lb}(\Hom_{gr}(\gmx{1},\gm)) \simeq \Z/p\Z$$ then
$\{\widetilde{\clE}^{(\mu,\lb;F)}\}=0$. \item[b)] Let us suppose
$v(\lb)=0$. Since $S_\lb=\emptyset$, then
$\{\widetilde{\clE}^{(\mu,\lb;F)}\}=0$.
\end{itemize}
\end{ex}

\subsection{Interpretation of
$\Ext^1(G_{\mu,1},\gm)$}

First of all, we   briefly recall a useful spectral sequence . Let
$\mathcal{E}xt^i(G,H)$ denote the fppf-sheaf on $Sch_{/S}$,
associated to the presheaf $X\longmapsto
\Ext^i{(G\times_{S}X,H{\times_S}X)}$. Then we have a spectral
sequence
$$
E_{2}^{ij}=H^i(S,\mathcal{E}xt^j(G,H))\Rightarrow \Ext^{i+j}(G,H),
$$
which  in low degrees gives
\begin{equation}\label{eq:spectral sequence per ext}
\begin{aligned}
0\too H^1(S,\mathcal{E}xt^0(G,H))&\too \Ext^{1}(G,H)\too
H^0(S,\mathcal{E}xt^1(G,H))\too\\
&\too H^2(S,\mathcal{E}xt^0(G,H))\too \Ext^{2}(G,H).
\end{aligned}
\end{equation}

All the groups of cohomology are calculated in the fppf topology.
Moreover $H^1(S,\mathcal{E}xt^0(G,H))$ is isomorphic to the
subgroup of $\Ext^1(G,H)$ formed by the extensions $E$ which split
over some faithfully flat affine $S$-scheme of finite type (cf.
\cite[III 6.3.6]{DG}). We suppose that  $G$ acts trivially on $H$,
then $\mathcal{E}xt^0(G,H)=\mathcal{H}om_{gr}(G,H)$. We will
consider the case $H=\gm$ and $G$ a finite flat group scheme. In
this case, $\mathcal{E}xt^0(G,\gm)$ is  by definition  the Cartier
dual of $G$, denoted by $G^\vee$. We recall the following result
which will play a role in the description of extensions of
$\gmx{1}$ by $\glx{1}$ (see \ref{teo:ext1(glx,gmx)} below).
\begin{thm}\label{prop:ext1(gmx,gm)}
  Let $G$ be a commutative finite flat group scheme over  $S$.
  Then
the canonical map
$$
H^1(S,G^\vee)\too \Ext^1(G,\gm)
$$
is bijective.
\end{thm}
\begin{proof}
This is a Theorem of S.U. Chase. For a proof see \cite{Wat2}. We
stress that he proves $\cal{E}xt^{1}(G,\gm)=0$, then he applies
\eqref{eq:spectral sequence per ext}. We remark  that he proves
everything in the fpqc site. However the same proof works in the
fppf site.
\end{proof}
We  apply this result to $G=\glx{1}$. We have that the map
\begin{equation*}\label{eq:ext1=H^1 per glbn}
H^{1}(S,\glx{1}^\vee){\too} \Ext^{1}({\glx{1}},{\gm})
\end{equation*}
is an isomorphism. Moreover the following is true.
\begin{prop}\label{prop:torsori schemi normali}
Let $X$ be a normal integral scheme. 
For any finite and flat commutative group scheme $G$ over $X$,
$$
i^*:\h{1}{X}{G}\too \h{1}{\Sp(K(X))}{G_{|K(X)}}
$$
is injective, where $i:\Sp(K(X))\too X$ is the generic point.
\end{prop}
\begin{proof} A sketch of the proof has been suggested to us by F.
Andreatta. We recall that for any commutative group scheme $G'$
over a scheme $X'$ we have that $\h{1}{X'}{G'}$ is a group and it
classifies the $G'$-torsors over $X'$. Suppose there exists a
$G$-torsor $f:Y\too X$ such that $i^*f:i^*Y\to \Sp(K(X))$ is
trivial. This means there exists a section $s$ of $i^{*}f$. We
consider the scheme $Y_0$ which is the closure of $s(\Sp(K(X)))$
in $Y$. Then $f_{|Y_0}:Y_0\too X$ is a finite birational morphism
with $X$ a normal integral scheme. So, by Zariski's Main Theorem
(\cite[4.4.6]{liu}), we have that $f_{|Y_0}$ is an open immersion,
and so it is an isomorphism. So we have a section of $f$ and $Y$
is a trivial $G$-torsor.
\end{proof}
\begin{rem}The hypothesis $G$ finite and $X$ normal   are necessary.
For the first it is sufficient to observe that  any $\gm$-torsor
is trivial on
$\Sp(K(X))$. 
For the second one, consider $X=\Sp(k[x,y]/(x^p-y^{p+1}))=\Sp(A)$,
with $k$ any field of characteristic $p>0$, and $Y$ the
$\alpha_p$-torsor $\Sp(A[T]/(T^{p}-y))$. Generically this torsor
is trivial since we have $y=(\frac{x}{y})^p$. But $Y$ is not
trivial since $y$ is not a $p$-power in $A$.
\end{rem}
\begin{cor}Let $X$ be a normal integral scheme.  Let $f:Y\too X$ be a
morphism with a rational section and let $g:G\too G'$ a map of
finite and flat commutative group schemes over $X$, which is an
isomorphism over $\Sp(K(X))$. Then
$$
f^*g_*:\hxg{G}\too H^1(Y,G'_{Y})
$$
is injective.
\end{cor}
\begin{proof}
By hypothesis $\Sp(K(X))\too X$ factorizes through $f:Y\too X$. If
$i:\Sp(K(X))\too X$,  we have
$$
i_*:\hxg{G}\on{g_*}{\too} \hxg{G'}\on{f^*}{\too} H^1(Y,G'_{Y})\too
H^1(\Sp(K(X)),G_{K(X)}).
$$
Therefore, by the previous proposition, it follows that
$$
\hxg{G}\too H^1(Y,G'_{Y})
$$
is injective.
\end{proof}
\begin{rem}The previous corollary can be applied, for instance, to the case
$f=\id_X$ or to the case $f:U\too X$ an open immersion and
$g=\id_G$. Roberts (\cite[p. 692]{Rob}) has proved the corollary
in the case $f=\id_X$, with $X=\Sp(A)$ and $A$ the integer ring of
a local number field.
\end{rem}

By \ref{prop:ext1(gmx,gm)} and  \ref{prop:torsori schemi normali}
we obtain the following result.
\begin{cor}\label{cor:restrizione di ext^1 iniettiva} Let $G$ be a commutative finite flat group scheme over  $S$. The restriction map
$$
\Ext^{1}(G,\gm)\too \Ext^{1}(G_{K},{\gm}_{|K})
$$
is injective.
\end{cor}

Let us consider a commutative finite and flat group scheme $G$ of
order $n$. We also consider  the $n^{th}$ power map $n:\gm\too
\gm$ . It induces a morphism $n_*:\Ext^{1}(G,{\gm})\too
\Ext^{1}({G},{\gm})$. We have the following commutative diagram
$$
\xymatrix{H^{1}(S,G^\vee)\ar[r]\ar[d]^{n_*} &
\Ext^{1}({G},{\gm})\ar[d]^{n_*}\\
H^{1}(S,G^{\vee})\ar[r]& \Ext^{1}(G,{\gm}), }
$$
where the horizontal maps are isomorphisms by
\ref{teo:ext1(glx,gmx)}. We remark that $n_*:H^{1}(S,G^\vee)\too
H^{1}(S,G^\vee)$ is the zero morphism since the  map
$n_*:G^\vee\too G^\vee$, induced by $n:\gm\too \gm$,  is the zero
morphism. This proves the following lemma.
\begin{lem}\label{lem:n*=0}
Let $G$ be a commutative finite and flat group scheme of order
$n$. Then
$$
n_*:\Ext^{1}(G,{\gm})\too \Ext^{1}({G},{\gm})
$$
is the zero morphism.
\end{lem}
\subsection{Description of $\Ext^1(\gmx{1},\glx{1})$}\label{par:Ext1(gmu,glb)}

We finally have all the ingredients to  give a description of the
group $\Ext^1(\gmx{1},\glx{1})$. In particular we will focus on
the extensions which are isomorphic, as group schemes, to
$\Z/p^2\Z$ on the generic fiber.

First of all we remark that if $v(\mu_1)=v(\mu_2)$ and
$v(\lb_1)=v(\lb_2)$ then
$$
\Ext^1(G_{\mu_1,1},G_{\lb_1,1})\simeq\Ext^1(G_{\mu_2,1},G_{\lb_2,1}).
$$
Indeed we know, by hypothesis, that  there exist two isomorphisms
\mbox{$\psi_1:G_{\lb_1,1}\too G_{\lb_2,1}$} and
$\psi_2:G_{\mu_2,1}\too G_{\mu_1,1}$.  Then we have that
$$
(\psi_1)_*\circ(\psi_2)^*:\Ext^1(G_{\mu_1,1},G_{\lb_1,1})\too\Ext^1(G_{\mu_2,1},G_{\lb_2,1})
$$
is an isomorphism.

We now recall what happens if $v(\mu)=v(\lb)=0$. In this case we
have the following result.
\begin{prop}\label{eq:ext1(mup,mup)}Let $A$ be a d.v.r or a field. Then there exists an exact sequence
$$
0\too\Z/p\Z\too \Ext^1_{A}(\mu_p,\mu_p)\too H^1(\Sp(A),\Z/p\Z)\too
0
$$
\end{prop}
\begin{proof}
The proposition  is proved in \cite[3.7]{SS8} when $A$ a is d.v.r.
The same proof works when $A$ is a field.
\end{proof}
%
Let us define the extension of $\mu_p$ by $\mu_p$
$$
\clE_{i,A}=\Sp(A[S_1,S_2]/(S_1^p-1,\frac{S_2^p}{S_1^i}-1)).
$$
It is the kernel of the morphism $(\gm)^2\too (\gm)^2$ given by
$(S_1,S_2)\too (S_1^p,S_1^{-1}S_2^p)$
Then  
the group $\Z/p\Z$ in the above proposition is formed by the
extensions $\clE_{i,A}$.
\begin{defn}\label{def:clE(mu,lb,F,i)}
 Let $F\in \Hom({\gmx{1}}_{|S_\lb},{\gm}_{|S_\lb})$, $j\in
\Z/p\Z$ such that  $$F(S)^p(1+\mu S)^{-j}=1\in
\Hom({\gmx{1}}_{|S_{\lb^p}},{\gm}_{|S_{{\lb^p}}}).$$ Let $\tilde{F}(S) \in R[S]$  be a  lifting of $F$. 
We denote by $\clE^{(\mu,\lb;F,j)}$  the subgroup scheme of
$\clE^{(\mu,\lb;F)}$ given
 on the level of schemes  by
$$
\clE^{(\mu,\lb;F,j)}=\Sp\bigg(R[S_1,S_2]/\big(\frac{(1+\mu
S_1)^p-1}{\mu^p},\frac{(\tilde{F}(S_1)+\lb S_2)^p(1+\mu
S_1)^{-j}-1}{\lb^p}\big)\bigg).
$$
\end{defn}

 We moreover  define the following
homomorphisms of group schemes
$$
\glx{1}\too \clE^{(\mu,\lb;F,j)}
$$
by
\begin{align*}
&S_1\longmapsto 0\\
&S_2\longmapsto S +\frac{1-\tilde{F}(0)}{\lb}
\end{align*}
and
$$
\clE^{(\mu,\lb;F,j)}\too G_{\mu,1}
$$
 by
\begin{align*}
&S\too S_1.
\end{align*}
It is easy to see that
$$
0\too \glx{1}\too \clE^{(\mu,\lb;F)}\too G_{\mu,1}\too 0
$$
is exact. A different choice of the lifting $\tilde{F}(S)$ gives
an isomorphic extension. It is easy to see that
$(\clE^{(\mu,\lb;F,j)})_K\simeq (\Z/p^2\Z)_K $, as a group scheme,
if  $j\neq 0$ and $(\clE^{(\mu,\lb;F,0)})_K\simeq (\Z/p\Z\times
\Z/p\Z)_K $.

%

\begin{rem}In the above definition the integer $j$ is uniquely
determined by $F\in \Hom({\gmx{1}}_{|S_\lb},{\gm}_{|S_\lb})$ if
and only if $\lb^p\nmid \mu$.
\end{rem}
\vspace{.5cm} From the exact sequence over $S_\lb$
$$ 0\too G_{\mu,1}\on{i}{\too}
\gmu\on{\psi_{\mu,1}}{\too}\gmup\too 0
$$
we have that
\begin{equation}\label{eq:ker i_*}
\ker \bigg(i^*:\Hom_{gr}(\gmu_{|S_\lb},{\gm}_{|S_\lb})\too
\Hom_{gr}({G_{\mu,1}}_{|S_\lb},{\gm}_{|S_\lb})\bigg)={\psi_{\mu,1}}_*\Hom_{gr}(\gmup_{|S_\lb},{\gm}_{|S_\lb})
\end{equation}
So let $F(S)\in
\Hom({\gmx{1}}_{|S_{\lb^p}},{\gm}_{|S_{{\lb^p}}})$. By
\ref{lem:suriettività mappa tra hom} we can choose a representant
of $F(S)$ in $\in \Hom({\gmu}_{|S_{\lb^p}},{\gm}_{|S_{{\lb^p}}})$
which we denote again $F(S)$ for simplicity. Therefore, by
\eqref{eq:ker i_*}, we have that $F(S)^p(1+\mu S)^{-j}=1\in
\Hom({\gmx{1}}_{|S_{\lb^p}},{\gm}_{|S_{{\lb^p}}})$ is equivalent
to saying that there
exists $G\in\Hom({\gmup}_{|S_{\lb^p}},{\gm}_{|S_{\lb^p}})$ with
the property that $F(S)^p(1+\mu S)^{-j}=G(\frac{(1+\mu
S)^p-1}{\mu^p})\in
\Hom({\gmu}_{|S_{\lb^p}},{\gm}_{|S_{{\lb^p}}})$. This implies that
$\clE^{(\mu,\lb;F,j)}$  can be seen as the kernel of the isogeny
\begin{align*}
\psi^j_{\mu,\lb,F,G}:\clE^{(\mu,\lb;F)}&\too \clE^{(\mu^p,\lb^p;G)}\\
          S_1&\longmapsto \frac{(1+\mu S_1)^p-1}{\mu^p}\\
          S_2&\longmapsto \frac{(\tilde{F}(S_1)+\lb
S_2)^p(1+\mu S_1)^{-j}-\tilde{G}(\frac{(1+\mu
S_1)^p-1}{\mu^p})}{\lb^p}
\end{align*}
where $\tilde{F},\tilde{G}\in R[T]$  are liftings of $F$ and $G$.

As remarked in \ref{rem:v(mu)ge v(lb)}, if $v(\mu)\ge v(\lb)$  we
can suppose
$$
\tilde{F}(S)=\sum_{i=0}^{p-1}\frac{a^i}{i!}S^i
$$
 with $a^p\equiv 0\mod \lb$.

 \begin{ex}\label{ex:Z/p^2Z SS}This example has been the main motivation for our definition of the group schemes $\clE^{(\mu,\lb;F,j)}$. Let us define
$$
\eta=\sum_{k=1}^{p-1}\frac{(-1)^{k-1}}{k}\lb_{(2)}^k.
$$
We remark that $v(\eta)=v(\lb_{(2)})$. We consider
$$
F(S)=E_p(\eta S)=\sum_{k=1}^{p-1}\frac{(\eta S)^k}{k!}
$$
It has been shown in \cite[\S 5]{SS4} that, using our notation,
$$
\Z/p^2\Z\simeq \clE^{(\lb_{(1)},\lb_{(1)};E_p(\eta S),1)}.
$$


A similar description of $\Z/p^2\Z$ was independently found by
Green and Matignon (\cite{GM1}).
\end{ex}
\vspace{.5cm}

\begin{ex}\label{ex:Gmun} It is easy to see that the group scheme $\glx{2}$ 
is isomorphic to
$$
\clE^{(\lb^p,\lb;1,1)}.
$$
Moreover if we have an extension  of type $\clE^{(\mu,\lb;F,j)}$
with $F(S)=1$ 
then  $v(\mu)\ge pv(\lb)$. Indeed we have that
$$
\clE^{(\mu,\lb;1,j)}=\Sp\bigg(A[S_1,S_2]/\bigg(\frac{(1+\mu
S_1)^p-1}{\mu^p},\frac{(1+\lb S_2)^p(1+\mu
S_1)^{-j}-1}{\lb^p}\bigg)\bigg).
$$
Since $(1+\lb S_2)^p=1\in
\Hom_{gr}({\gmx{1}}_{|S_{\lb^p}},{\gm}_{|S_\lb^p})$ then $${(1+\lb
S_2)^p(1+\mu
S_1)^{-j}}=1\in\Hom_{gr}({\gmx{1}}_{|S_{\lb^p}},{\gm}_{|S_{\lb^p}})$$
if and only if $v(\mu)\ge pv(\lb)$. In particular we remark that,
in such a  case, $v(\lb)\le v(\lb_{(2)})$. Otherwise
$$pv(\lb)>pv(\lb_{(2)})=v(\lb_{(1)})\ge v(\mu),$$ which is not
possible.
\end{ex}
\vspace{.3cm}

 Let us define, for any $\mu,\lb\in R$ with $v(\mu),v(\lb)\le v(\lb_{(1)})$, the group
\begin{align*}
rad_{p,\lb}(<1+\mu &S>):=\bigg\{(F(S),j)\in
\Hom_{gr}({\gmx{1}}_{|S_\lb},{\gm}_{|S_\lb})\times \Z/p\Z \text{
such that }\\
& F(S)^p(1+\mu S)^{-j}=1\in
\Hom({\gmx{1}}_{|S_{\lb^p}},{\gm}_{|S_{{\lb^p}}})\bigg\}/<(1+\mu
S,0)>.
\end{align*}

We define $$ \beta: rad_{p,\lb}(<1+\mu S>) \on{}{\too}
\Ext^{1}(G_{\mu,1},\glx{1})$$ by
$$ (F(
S),j)\longmapsto \clE^{(\mu,\lb;F( S),j)}
$$
We remark that the image of $\beta$ is the set
$\{\clE^{(\mu,\lb;F( S),j)}\}$.

\begin{lem}
$\beta$ is a morphism of groups. In particular the set
$\{\clE^{(\mu,\lb;F(S),j)}\}$ is a subgroup of
$\Ext^{1}(G_{\mu,1},\glx{1})$.
\end{lem}
\begin{proof}
Let $i:G_{\lb,1}\too \glb$. We remark that $$ i_*(\beta(F,j))=i_*(
\clE^{(\mu,\lb;F( S),j)})=\tilde{\clE}^{(\mu,\lb;F)}=\delta(F)$$
for any $(F,j)\in rad_{p,\lb}(<1+\mu S>)$. Moreover  by
construction
$$
(\clE^{(\mu,\lb;F( S),j)})_K=(\cal{E}_{j})_K\in
\Ext^{1}({\mu_p}_{,K},{\mu_p}_{,K}).
$$

 Let $(F_1,j_1),(F_2,j_2)\in rad_{p,\lb}(<1+\mu S>)$. Then
\begin{equation}\label{eq:i_*=0}
i_*(\beta(F_1,j_1)+\beta(F_2,j_2)-\beta(F_1+F_2,j_1+j_2))=\delta(F_1)+\delta(F_2)-\delta(F_1+F_2)=0
\end{equation}
since $\delta$ is a morphism of groups. And
\begin{equation}\label{eq:restrizione su K e' morfismo}
(\beta(F_1,j_1)+\beta(F_2,j_2)-\beta(F_1+F_2,j_1+j_2))_K=\clE_{j_1,K}+\clE_{j_2,K}-\clE_{j_1+j_2,K}=0,
\end{equation}
since  $\Z/p\Z\simeq\Ext^1(\mu_{p,K},\mu_{p,K})$ through the map
$j\mapsto \clE_{j,K}$.
 By
\eqref{eq:i_*=0} it follows that
$$
\beta(F_1,j_1)+\beta(F_2,j_2)-\beta(F_1+F_2,j_1+j_2)\in\ker i_*.
$$
and then,  by \eqref{eq:ker phi} and \eqref{eq:def delta'}, we
have
$$
\beta(F_1,j_1)+\beta(F_2,j_2)-\beta(F_1+F_2,j_1+j_2)=(\sigma_j)^*\Lambda.
$$
for some $j\in \Z/p\Z$. By \eqref{eq:restrizione su K e' morfismo}
it follows that
$$
((\sigma_j)^*\Lambda)_K=\clE_{j,K}=0,
$$
therefore $j=0$. So $\beta$ is a morphism of groups.
The last assertion is clear.
\end{proof}
 We  now give a description of $\Ext^{1}(G_{\mu,1},\glx{1})$.
\begin{thm}\label{teo:ext1(glx,gmx)}Suppose that $\lb,\mu\in R$ with $v(\lb_{(1)})\ge v(\lb),v(\mu)$.
 The 
following sequence
\begin{equation*}
\begin{array}{ll}
0\too rad_{p,\lb}(<1+\mu S>) \on{\beta}{\too}&
\Ext^{1}(G_{\mu,1},\glx{1})\on{\alpha^{\lb}_*\circ i_*}\too\\
&\too \ker \bigg(H^1(S,G_{\mu,1}^\vee)\too
H^1(S_\lb,G_{\mu,1}^\vee)\bigg)
\end{array}
\end{equation*}
is exact. In particular $\beta$ induces an isomorphism
$rad_{p,\lb}(<1+\mu S>)\simeq\{\clE^{(\mu,\lb;F,j)}\}$.
\end{thm}
\begin{proof}
Using \eqref{eq:ker alpha} and  \ref{prop:ker alpha=clE}, we consider the following commutative diagram 
\begin{equation}\label{eq:succ esatta per ext1} \xymatrix{0\ar[r]
&\{\widetilde{\clE}^{(\mu,\lb;F)}\}\ar[r]\ar^{\widetilde{\psi_{{\lb,1}_*}}}[d]&
\Ext^{1}(\gmx{1},\glb)\ar^{\alpha_*^{\lb}}[r]\ar^{\psi_{{\lb,1}_*}}[d]&\Ext^{1}(\gmx{1},\gm)\ar^{p_*}[d]\ar[r]&  \Ext^{1}({G_{\mu,1}}_{|S_\lb},{\gm}_{|S_\lb})\ar^{p_*}[d]\\
 0\ar[r]  &\{\widetilde{\clE}^{(\mu,\lb^p;G)}\}\ar[r]&
\Ext^{1}(\gmx{1},\g^{(\lb^p)})\ar^{{\alpha_*^{\lb^p}}}[r]
&\Ext^{1}(\gmx{1},\gm)\ar[r]      &
\Ext^{1}({G_{\mu,1}}_{|S_\lb},{\gm}_{|S_\lb})}
\end{equation}
 The map
$\widetilde{\psi_{{\lb,1}_*}}$, 
induced by ${\psi_{\lb,1}}_*:\Ext^{1}(\gmx{1},\glb)\too
\Ext^{1}(\gmx{1},\glbp)$, is given by
$\widetilde{\clE}^{(\mu,\lb;F)}\longmapsto
\widetilde{\clE}^{(\mu,\lb^p;F^p)}$. 
Now, since  $\gmx{1}$ is of order $p$ then,
$p_*:\Ext^1(\gmx{1},\gm)\to
\Ext^1(\gmx{1},\gm)$ 
is the zero map (see \ref{lem:n*=0}).
Moreover, by \eqref{eq:spectral sequence per ext} and
\ref{prop:ext1(gmx,gm)}, we have the following situation
$$
\xymatrix@1{0\ar[r]&H^1(S,G_{\mu,1}^\vee)\ar[r] \ar[d]& \Ext^1(\gmx{1},\gm) \ar[d]^{}\ar[r]&0\\
0\ar[r]&H^1(S_\lb,G_{\mu,1}^\vee)\ar[r]&\Ext^1({\gmx{1}}_{|S_\lb},{\gm}_{|S_\lb})}
$$
which
 implies that
$ Im(\alpha^\lb_*)\simeq \ker
(H^1(S,G_{\mu,1}^\vee)\too H^1(S_\lb,G_{\mu,1}^\vee))$. 
%

So applying the 
snake lemma to  \eqref{eq:succ esatta per ext1} we obtain 
\begin{equation}\label{eq:snake}
\begin{array}{l}
0\too \ker(\widetilde{\psi_{{\lb,1}_*}})
\on{\widetilde{\delta}}{\too}
\ker(\psi_{{\lb,1}_*})\on{\alpha^{\lb}_*} \too \ker
\bigg(H^1(S,G_{\mu,1}^\vee)\too
H^1(S_\lb,G_{\mu,1}^\vee)\bigg). 
\end{array}
\end{equation}
We now divide the proof in some steps.

\begin{tabular}{|c|}
  \hline
  Connection between $\ker(\widetilde{{\psi_{\lb,1}}_*})$ and $rad_{p,\lb}(<1+\mu S>)$.\\
  \hline
\end{tabular}
We are going to give the connection in the form of the isomorphism
\eqref{eq:isomorfismo} below. We recall that, by \eqref{eq:ker
phi}, $i:\glx{1}\too \glb$ induces an isomorphism
\begin{equation}\label{eq:ker phi bis}
i_*:\Ext^{1}(G_{\mu,1},\glx{1})/\delta'(\Hom_{gr}(G_{\mu,1},\glbp))\too
\ker({\psi_{\lb,1}}_*);
\end{equation}
for the definition of $\delta'$ see \eqref{eq:def delta'}.

By  \ref{prop:ker alpha=clE} we have an isomorphism 
\begin{align*}
    \delta:\Hom_{gr}({\gmx{1}}_{|S_{\lb}},{\gm}_{|S_\lb})/r_{\lb'}(\Hom_{gr}(\gmx{1},\gm))\too \{\widetilde{\clE}^{(\mu,\lb;F)} \}
\end{align*}
Through this identification we can identify 
$\ker(\widetilde{{\psi_{\lb,1}}_*})$ with
\begin{equation} \label{eq:ker(tilde psi)}
\begin{aligned}
\bigg\{F(S)\in \Hom_{gr}({\gmx{1}}_{|S_\lb},{\gm}_{|S_\lb})|
\exists &i\in r_{\lb^p}(\Hom_{gr}(\gmx{1},\gm)) \text{ such that }\\
 F(S)^p(1+\mu S)^{-i}=1\in
\Hom({\gmx{1}}_{|S_{\lb^p}}&,{\gm}_{|S_{{\lb^p}}})\bigg\}/<1+\mu
S>.
\end{aligned}
\end{equation}
Moreover 
\begin{equation}\label{eq:inj}
\widetilde{\delta}:\ker(\widetilde{{\psi_{\lb,1}}_*})\hookrightarrow
\Ext^{1}(G_{\mu,1},\glx{1})/\delta'(\Hom_{gr}(G_{\mu,1},\glbp))\In
\Ext^{1}(\gmx{1},\g^{(\lb)})
\end{equation}
 is defined by 
$\widetilde{\delta}(F)=\delta(F)=\widetilde{\clE}^{(\mu,\lb;F)}$. 

We now define  a morphism of groups
$$\iota:\ker(\widetilde{{\psi_{\lb,1}}_*})\too
r_{\lb^p}(\Hom_{gr}(\gmx{1},\gm))$$ as follows: for  any $F(S)\in
\ker(\widetilde{{\psi_{\lb,1}}_*})$,  $\iota(F)=i_F$  is the
unique $i\in r_{\lb^p}(\Hom_{gr}(\gmx{1},\gm))$ such that
$F(S)^p(1+\mu S)^{-i}=1\in
\Hom({\gmx{1}}_{|S_{\lb^p}},{\gm}_{|S_{{\lb^p}}})$. The morphism
of groups
\begin{align}\label{eq:isomorfismo}
\ker(\widetilde{\psi_{{\lb,1}_*}})\times \Hom_{gr}
(G_{\mu,1},\glbp)&\too rad_{p,\lb}(<1+\mu S>)\\
(F,j)&\longmapsto (F,i_F+j)\nonumber
\end{align}
is  an isomorphism. We prove only the surjectivity  since the
injectivity is clear. Now, if $\lb^p\nmid \mu$ then
$\Hom_{gr}(\gmx{1},\glbp)=0$ and
$r_{\lb^p}(\Hom_{gr}(\gmx{1},\gm))=\Z/p\Z$. So, if $(F,j)\in
rad_{p,\lb}(<1+\mu S>)$, then $j\in
r_{\lb^p}(\Hom_{gr}(\gmx{1},\gm))$. So $i_F=j$. Hence
$(F,0)\mapsto (F,i_F)=(F,j)$. While if $\lb^p\mid \mu$ then
$\Hom_{gr}(\gmx{1},\glbp)=\Z/p\Z$ and
$r_{\lb^p}(\Hom_{gr}(\gmx{1},\gm))=0$. Hence
$$\ker(\widetilde{{\psi_{\lb,1}}}_*)=\bigg\{F(S)\in
\Hom_{gr}({\gmx{1}}_{|S_\lb},{\gm}_{|S_\lb})|
 F(S)^p=1\in
\Hom({\gmx{1}}_{|S_{\lb^p}},{\gm}_{|S_{{\lb^p}}})\bigg\}.$$ Let us
now take $(F,j)\in rad_{p,\lb}(<1+\mu S>).$ This means that
$$F(S)^p=(1+\mu S)^j=1\in \Hom({\gmx{1}}_{|S_{\lb^p}},{\gm}_{|S_{{\lb^p}}}).$$
Therefore $F(S)\in \ker(\widetilde{\psi_{{\lb,1}_*}})$ and
$i_F=0$. So
$$
(F,j)\longmapsto (F,i_F+j)=(F,j).
$$

\begin{tabular}{|c|}
  \hline
  Interpretation of $\beta$.\\
  \hline
\end{tabular}
We now define the morphism of groups
\begin{align*}
\varrho:\ker(\widetilde{{\psi_{\lb,1}}_*})&\too
\Ext^{1}(G_{\mu,1},\glx{1})\\
        F&\longmapsto\beta(F,i_F)= \clE^{(\mu,\lb;F,i_F)}
\end{align*}
We recall the definition of $\delta'$ given in \eqref{eq:def
delta'}: $$\delta':\Hom_{gr} (G_{\mu,1},\glbp)\too
\Ext^{1}(G_{\mu,1},\glx{1})$$ is defined by
$\delta'(\sigma_i)=\sigma_i^*(\Lambda)$.
 Then, under the
isomorphism \eqref{eq:isomorfismo}, we have
$$
\beta=\rho+\delta':\ker(\widetilde{\psi_{{\lb,1}_*}})\times
\Hom_{gr} (G_{\mu,1},\glbp) {\too} \Ext^{1}(G_{\mu,1},\glx{1})
$$
%

%
\begin{tabular}{|c|}
  \hline
   Injectivity of $\beta$.\\
  \hline
\end{tabular}
First of all we observe that $\widetilde{\delta}$ factors through
$\rho$, i.e.
\begin{equation}\label{eq:delta=i circ rho}
\widetilde{\delta}=i_*\circ
\rho:\ker(\widetilde{{\psi_{\lb,1}}_*})\on{\rho}{\too}
\Ext^{1}(G_{\mu,1},\glx{1})\on{i_*}{\too}
\ker({\psi_{{\lb,1}_*}}).
\end{equation}
Indeed 
$$
i_*\circ
\rho(F)=i_*(\clE^{(\mu,\lb;F,i_F)})=\widetilde{\clE}^{(\mu,\lb,F)}=\widetilde{\delta}(F).
$$ 
In particular, since $\tilde{\delta}$ is injective, $\rho$ is
injective, too.

We now prove that $\beta=\rho+\delta'$ is injective, too. By
\eqref{eq:ker
phi bis}, 
$$i_*\circ \delta'=0.$$ Now, if
$(\rho+\delta')(F,\sigma_i)=0$, then $\rho(F)=-\delta'(\sigma_i)$.
So
$$
\widetilde{\delta}(F)=i_*(\rho(F))=i_*(-\delta'(\sigma_i))=0.
$$
But $\widetilde{\delta}$ is injective, so $F=1$. Hence
$\delta'(\sigma_i)=0$. But by \eqref{eq:ker phi}, also $\delta'$
is injective. Then $\sigma_i=0$.

\begin{tabular}{|c|}
  \hline
 Calculation of  $Im \beta$.\\
  \hline
\end{tabular}
 We finally prove  $Im(\rho+\delta')=  \ker({\alpha^{\lb}_*}\circ
i_*)$. Since $\widetilde{\delta}=i_*\circ \rho$,
$\alpha^{\lb}_*\circ\widetilde{\delta}=0$ and $i_*\circ\delta'=0$
then
$$
\alpha^{\lb}_*\circ i_*\circ(\rho+\delta')=\alpha^\lb_*\circ
i_*\circ\rho+\alpha^{\lb}_*\circ(
i_*\circ\delta')=\alpha^{\lb}_*\circ
i_*\circ\rho=\alpha^{\lb}_*\circ\widetilde{\delta}=0.
$$
So $Im(\rho+\delta')\In  \ker({\alpha^{\lb}_*}\circ i_*)$.  On the
other hand, if  $E\in \Ext^{1}(G_{\mu,1},\glx{1})$ is such that
${\alpha^{\lb}_*}\circ i_*(E)=0$, then, by \eqref{eq:snake}, there
exists $F\in \ker(\widetilde{\psi_{{\lb,1}_*}})$ such that
$i_*(E)=\widetilde{\delta}(F)=i_*(\rho(F))$. Hence, by
\eqref{eq:ker phi bis}, $E-\rho(F)\in Im(\delta')$. Therefore
$Im(\rho+\delta')= \ker({\alpha^{\lb}_*}\circ i^*)$. Moreover
since $i_*:\Ext^{1}(G_{\mu,1},\glx{1})\too \ker({\psi_{\lb,1}}_*)$
is surjective then $Im(\alpha^{\lb}_*)=Im(\alpha^{\lb}_*\circ
i_*)$. We have so proved, using also \eqref{eq:snake}, that the
following sequence
\begin{align*}
0\too \ker(\widetilde{\psi_{{\lb,1}_*}})\times \Hom_{gr}
(G_{\mu,1},\glbp)& \on{\rho+\delta'}{\too}
\Ext^{1}(G_{\mu,1},\glx{1})\on{\alpha^{\lb}_*\circ i_*}\too\\
\too &\ker \bigg(H^1(S,G_{\mu,1}^\vee)\too
H^1(S_\lb,G_{\mu,1}^\vee)\bigg)
\end{align*}
is exact.
Finally, by definitions, it follows that
$$
\beta(rad_{p,\lb}(<1+\mu S>))=\{\clE^{(\mu,\lb;F,j)}\}.
$$
\end{proof}
\begin{ex}\label{ex:teorema per v(lb)=0}Let us suppose $v(\lb)=0$. In such a case  $rad_{p,\lb}(<1+\mu T>)=\Z/p\Z$.
Hence by the theorem we have
$$
0\too \{\clE^{(\mu,\lb;1,j)}|j\in \Z/p\Z\}\too
\Ext^{1}(G_{\mu,1},\mup)\too H^1(S,G_{\mu,1}^{\vee})\too 0.
$$


\end{ex}
\begin{ex}Let us now suppose $v(\mu)=0$ and $v(\lb)>0$. In such a case
$$\Hom_{gr}({\mup}_{|S_\lb},{\gm}_{|S_\lb})=<1+\mu T>.$$ Hence it is
easy to see that
$$
rad_{p,\lb}(<1+\mu T>)=0.
$$ Therefore, by the
theorem,
$$
\Ext^{1}(\mup,G_{\lb,1})\too \ker \bigg(H^1(S,\Z/p\Z)\too
H^1(S_\lb,\Z/p\Z)\bigg)
$$
is an isomorphism.
\end{ex}
\begin{cor}\label{cor:ext1(gmx,glx) a meno di estendere R}
Under the hypothesis of the theorem, any extension $E\in
\ext^1(\gmx{1},\glx{1})$ is of type
$\clE^{(\mu,\lb;F,j)}$, up to an extension of $R$. 
In particular any extension is commutative.
\end{cor}
\begin{proof}
Let $E\in \ext^1{(\gmx{1},\glx{1})}$. Suppose  that
$\alpha^{\lb}_*(i_*E)=[S']$, with $S'\too S$ a
$G_{\mu,1}^\vee$-torsor.  
We consider the integral closure  $S''$ of $S$ in $S'_K$. Up to a localization  (in the case $S''\too S$ is étale), we can suppose $S''$ local. 
So $S''=\Sp(R'')$ where $R''$ is a no\oe therian local integrally
closed ring of dimension $1$, i.e. a d.v.r. (see \cite[9.2]{mac}).
Since $S''_K\simeq S'_K$, then  $S'_K\times_K S''_K$ is a trivial
$E_K$-torsor over $S''_K$. By \ref{prop:torsori schemi normali} we
have that  $S'\times_S S''$ is a trivial $E$-torsor trivial over
$S''$. 
So, if we make the base change $f:S''\to S$, then
$\alpha^{\lb}_*(i_*(E_{S'}))=0$. By \ref{teo:ext1(glx,gmx)}, this
implies that
 $E''$ is of type
$$
\clE^{(\mu,\lb,F,j)}.
$$
 Hence any $E\in
\Ext^{1}(\gmx{1},\glx{1})$ is a commutative group scheme over an
extension $R'$ of $R$. So it is a commutative group scheme over
$R$.

\end{proof}

\begin{rem}\label{rem:essenzialmente tutti i gruppi di
ordine p2} Any $R$-group scheme  of order $p^2$ is of type $
\clE^{(\mu,\lb,F,j)}$, possibly after an extension of $R$. Indeed,
up to an extension of $R$, the generic fiber of any $R$-group scheme
is a constant group. Then the thesis follows from \ref{lem:modelli
di Z/p^2 Z sono estensioni} and \ref{cor:ext1(gmx,glx) a meno di
estendere R}.
\end{rem}

By \ref{prop:hochschild and extensions}  the extensions of
$\Z/p\Z$ by $\Z/p \Z$ over $K$, which are extensions of abstract
groups, are classified by $H_0^2(\Z/p\Z,\Z/p\Z)\simeq \Z/p\Z$ (see
for instance \cite[2.7]{SS8}). This group is formed by
$\clE_{j,K}$ with $j\in \Z/p\Z$. If $j\neq
0$ we have that $\clE_{j,K}$ 
is isomorphic, as a  group scheme, to $\Z/p^2\Z$, while if $j=0$ it is isomorphic to $  \Z/p\Z\times \Z/p\Z$. 
We also define    the following morphism of extensions
\begin{equation}
\begin{aligned}\label{eq:morfismo Z/p2Z to E}
\alpha_{\mu,\lb}:\clE^{(\mu,\lb;F,j)}&\too \clE_{j,R} \\
  S_1&\longmapsto 1+\mu S_1 \\
  S_2&\longmapsto F(S_1)+\lb S_2. \\ 
\end{aligned}
\end{equation}
It is an isomorphism on the generic fiber.
Now, by the theorem, we get that $\clE^{(\mu,\lb;F,j)}$ are the
only extensions which are isomorphic to $\clE_{j,K}$ on the
generic fiber.
\begin{cor} \label{cor:modelli di Z/p^2Z} 
The  extensions of type $\clE^{(\mu,\lb;F,j)}$ are the only
extensions $\clE\in \ext^{1}(\gmx{1},\glx{1})$
 which are isomorphic, as extensions, to $\clE_{j,K}$ on the generic fiber. In
 particular they are the unique finite and flat $R$-group schemes of order $p^2$ which are models of constant groups. More precisely, they are
 isomorphic on the generic fiber, as group schemes,
 to $\Z/p^2\Z$ if $j\neq 0$ and to $\Z/p\Z\times \Z/p\Z$ if $j=0$.
\end{cor}
\begin{proof}
As remarked above any $\clE^{(\mu,\lb;F,j)}$
has  the properties of the statement. We now prove that they are
the unique extensions of $G_{\mu,1}$ by $G_{\lb,1}$ to have these
properties.
 Let $\clE\in \ext^{1}(\gmx{1},\glx{1})$ be such that  $\clE_K\simeq \clE_{j,K}$ as group schemes. By \ref{prop:ext1(gmx,gm)}, \ref{eq:ext1(mup,mup)}  and \ref{teo:ext1(glx,gmx)}  we
have the following commutative diagram
$$
\xymatrix{\Ext^{1}(G_{\mu,1},\glx{1})\ar[r]^(0.36){\alpha^{\lb}_*\circ
i_*}\ar[d]& \ker \bigg(H^1(S,G_{\mu,1}^\vee)\too
H^1(S_\lb,G_{\mu,1}^\vee)\bigg)\ar[d]\\
\ext^1_{K}(\mu_p,\mu_p)\ar[r]^(.3){\alpha^{\lb}_*\circ i_*}&
\ext^1_{K}(\mu_p,\gm)\simeq H^1(\Sp(K),\Z/p\Z)\ar[r]& 0}
$$
where the vertical maps are the restrictions to the generic fiber.
Suppose now that $\clE_K$ is of type $\clE_{j,K}$. By
\ref{eq:ext1(mup,mup)} it follows that $\alpha^{\lb}_*\circ
i_*(\clE_K)=0$. Since the above diagram commutes, this means that
$({\alpha^{\lb}_*\circ i_*(\clE)})_K=0$. 
By \ref{prop:torsori schemi normali} we have that the second
vertical map of the diagram is injective. This means that
$$
{\alpha^{\lb}_*\circ i_*(\clE)}=0.
$$
So \ref{teo:ext1(glx,gmx)} 
implies that $\clE$ is of type $\clE^{(\mu,\lb;F,j)}$. Now, if $G$
is a model of a constant group, by \ref{lem:modelli di Z/p^2 Z
sono estensioni} we have that $G$ is an extension $\clE$ of
$G_{\mu,1}$ by $G_{\lb,1}$. Moreover, since $\clE_K$  is a
constant group, then $\clE_K\in \Ext^1(\Z/p\Z,\Z/p\Z)$. Therefore
$\clE_K\simeq \clE_j$
 for some $j$. So, by what we just proved, $\clE$ is of type $\clE^{(\mu,\lb;F,j)}$. The last assertion is clear.
\end{proof}
\subsection{$\Ext^1(G_{\mu,1},\glx{1})$ and the Sekiguchi-Suwa
theory}

 We now give a  description of $\clE^{(\mu,\lb;F,j)}$
through the Sekiguchi-Suwa theory. We study separately the cases
$\lb \nmid \mu$ and $\lb|\mu$.
\begin{cor}\label{prop:rad_p} 
Let $\mu ,\lb \in R$ be with $v(\lb_{(1)})\ge v(\lb)>v(\mu)$.
Then, no $\clE\in \Ext^1(G_{\mu,1},G_{\lb,1})$ is  a model of
$(\Z/p^2\Z)_K$. Moreover, if $p>2$ and $v(\mu)>0 $, the group
$\{\clE^{(\mu,\lb;F,j)} \} $ is isomorphic to 
\begin{align*} \bigg\{\bbf{a}\in
{\widehat{W}(R/\lb R)}^{\fr-[{\mu^{p-1}}]}|\exists
 \bbf{b} \in
{\widehat{W}(R/\lb^p R)}^{\fr-[{\mu^{p(p-1)}}]} \text{ such}
\text{
that }&\\
p\bbf{a}=[\frac{p}{\mu^{p-1}}]\textit{\textbf{b}}+V(\textit{\textbf{b}})\in
{\widehat{W}(R/\lb^p R)} \bigg\}\bigg/
 <[\mu],\bigg\{ [\frac{p}{\mu^{p-1}}]\textit{\textbf{b}}+
V(\textit{\textbf{b}}) |\bbf{b} \in& {\widehat{W}(R/\lb^p
R)}^{\fr-[{\mu^{p(p-1)}}]}\bigg\}>,
\end{align*}
through the map
$$
\bbf{a}\longmapsto \clE^{(\mu,\lb;E_p(\bbf{a},\mu; S),0)}.
$$
\end{cor}
\begin{rem}We know by \ref{pro:Hom(gmu,gm)}, \ref{lem:suriettività mappa tra hom} and  \ref{es:esempi di Ep(a,mu,T)} that any element of the set
defined above can be chosen of the type $[a]$ for some $a\in
(R/\lb R)^{\fr-[{\mu^{p-1}}]}$. So, if we have two elements as
above of the form $[a]$ and $[b]$ then $[a]+[b]=[c]$ for  some
$c\in (R/\lb R)^{\fr-[{\mu^{p-1}}]}$. We are not able to describe
explicitly this element. If we were able to do it we could have a
simpler description of the above set, as it happens in the case
$v(\mu)\ge v(\lb)$. We will see this in \ref{cor:clE se lb divide
mu}.

\end{rem}
\begin{proof}
We now prove the first statement. We remark that by
\ref{cor:modelli di Z/p^2Z} it is sufficient to prove the
statement only for the extensions in $\{\clE^{(\mu,\lb;F,j)}\}$.
 Let us consider the restriction map
$$
r:\Hom_{gr}({\gmx{1}}_{|S_{\lb^p}},{\gm}_{|S_{\lb^p}})\too
\Hom_{gr}({\gmx{1}}_{|S_{\lb}},{\gm}_{|S_{\lb}}).
$$
 The morphism $p:\Hom_{gr}(\gmu_{|S_{\lb}},{\gm}_{|S_\lb}){\too}\Hom_{gr}({\gmu}_{|S_{\lb^p}},{\gm}_{|S_{\lb^p}})$
defined  in \eqref{eq:def p} is given by $F(S)\mapsto F(S)^p$ and
induces a map
$$\Hom_{gr}({\gmx{1}}_{|S_{\lb}},{\gm}_{|S_\lb})\on{p}{\too}\Hom_{gr}({\gmx{1}}_{|S_{\lb^p}},{\gm}_{|S_{\lb^p}}).$$
Then
$$
\Hom_{gr}({\gmx{1}}_{|S_{\lb}},{\gm}_{|S_\lb})\on{p}{\too}\Hom_{gr}({\gmx{1}}_{|S_{\lb^p}},{\gm}_{|S_{\lb^p}})\on{r}{\too}
\Hom_{gr}({\gmx{1}}_{|S_{\lb}},{\gm}_{|S_{\lb}})
$$
 is the trivial morphism. Indeed
$$
(r\circ p)(F(S))=F(S)^p\in
\Hom_{gr}({\gmx{1}}_{|S_{\lb}},{\gm}_{|S_{\lb}}),
$$
which is zero by definition of group scheme morphisms and by the
fact that $G_{\mu,1}$ has order $p$. Now let us take $$F(S)\in
rad_{p,\lb}(<1+\mu S>)\simeq \{\clE^{(\mu,\lb;F,j)}\}.$$ By
definition
$$
F(S)^p(1+\mu S)^{-j}= 1 \in
\Hom({\gmx{1}}_{|S_{\lb^p}},{\gm}_{|S_{{\lb^p}}}),
$$
for some $j\in \Z/p\Z$. Hence
$$
r(F(S)^p(1+\mu S)^{-j})=(1+\mu S)^{-j}=1\in
\Hom_{gr}({\gmx{1}}_{|S_{\lb}},{\gm}_{|S_{\lb}}).
$$
If $(\clE^{(\mu,\lb;F,j)})_K\simeq (\Z/p^2\Z)_K$ then $j\neq 0$.
Therefore
$$
(1+\mu S)^{-j}=1\in
\Hom_{gr}({\gmx{1}}_{|S_{\lb}},{\gm}_{|S_{\lb}})
$$
 means $v(\mu)\ge v(\lb)$. So, if $v(\mu)<v(\lb)$, necessarily
 $j=0$. Hence
\begin{align*}
rad_{p,\lb}(<1+\mu S>):=\bigg\{&F(S)\in
\Hom_{gr}({\gmx{1}}_{|S_\lb},{\gm}_{|S_\lb}) \text{
such that }\\
& F(S)^p=1\in
\Hom({\gmx{1}}_{|S_{\lb^p}},{\gm}_{|S_{{\lb^p}}})\bigg\}/<1+\mu
S>.
\end{align*}
Therefore by \ref{teo:ext1(glx,gmx)}, \ref{lem:suriettività mappa
tra hom} and \eqref{eq:p a} we have the thesis.
\end{proof}

\begin{cor}\label{cor:clE se lb divide mu}
Let us suppose $p>2$. Let $\mu, \lb \in R\setminus \{0\}$ be with
$v(\lb_{(1)})\ge v(\mu)\ge v(\lb)$. Then, $ \{\clE^{(\mu,\lb;F,j)}
\} $ is isomorphic to the group
\begin{align*}
\Phi_{\mu,\lb}:=\bigg\{(a,j)\in (R/\lb R)^{\fr}\times\Z/p\Z&
\text{ such that } pa-j\mu=\frac{p}{\mu^{p-1}}a^p\in R/\lb^p
R\bigg\},
\end{align*}
through the map
$$
(a,j)\longmapsto
\clE^{(\mu,\lb;\sum_{i=0}^{p-1}\frac{a^i}{i!}S^i,j)}.
$$
\end{cor}
\begin{rem}\label{rem:0 in Phi}
It is clear that if $(0,j)\in \Phi_{\mu,\lb}$, with $j\neq 0$,
then $\mu\equiv 0 \mod \lb^p$.
\end{rem}

\begin{proof}
By \ref{lem:suriettività mappa tra hom}, \ref{cor: hom gmu gm se
lb divide mu}, \eqref{eq:p a} and \ref{ex:teorema per v(lb)=0}
(for the case $v(\lb)=0$) it follows that $rad_{p,\lb}(<1+\mu S>)$
is isomorphic to
\begin{equation}
\begin{aligned}\label{eq:rad p} \bigg\{(a,j)\in
(R/\lb R)^{\fr}\times \Z/p\Z|\exists  \bbf{b} \in
{\widehat{W}(R/\lb^p R)}^{\fr}
\\ \text{ such} \text{ that }
p[a]-j[\mu]=&[\frac{p}{\mu^{p-1}}]\textit{\textbf{b}}+V(\textit{\textbf{b}})\in
{\widehat{W}(R/\lb^p R)} \bigg\}.
\end{aligned}
\end{equation}
Let $a,j$ and $\textit{\textbf{b}}=(b_0,b_1,\dots)$ be as above.

By \cite[5.10]{SS4},
\begin{equation*}
p[a]\equiv (pa,a^p,0,\dots)\mod p^2.
\end{equation*}

Since $[\mu]\in \widehat{W}(R/\lb^p R)^{\fr}$ it follows by
\ref{lem:somma termine per termine} that
\begin{equation}\label{eq:i[mu]=[i mu]} j[\mu]= [j\mu]
\end{equation}
and \begin{equation}\label{eq:pa-kmu}
p[a]-j[\mu]=(pa-j\mu,a^p,0,0,\dots,0,\dots)\in
\widehat{W}(R/\lb^pR).
\end{equation}
We recall that$$
[\frac{p}{\mu^{p-1}}]\textit{\textbf{b}}=(\frac{p}{\mu^{p-1}}b_0,
\dots, (\frac{p}{\mu^{p-1}})^{p^i} b_i,\dots),$$ then, again by
\ref{lem:somma termine per termine}, we have
\begin{equation}\label{eq:p mu^(p-1) b=...b_0}
[\frac{p}{\mu^{p-1}}]\textit{\textbf{b}}+V(\textit{\textbf{b}})=
(\frac{p}{\mu^{p-1}} b_0, (\frac{p}{\mu^{p-1}})^pb_1+
b_0,\dots,(\frac{p}{\mu^{p-1}})^{p^{i+1}}b_{i+1}+b_i,\dots).
\end{equation}


Since $\bbf{b} \in {\widehat{W}(R/\lb^p R)} $ there exists $r\ge
0$ such that $b_i=0$ for any $i\ge r$. Moreover, comparing
\eqref{eq:pa-kmu} and \eqref{eq:p mu^(p-1) b=...b_0}, it follows
\begin{align*}
&(\frac{p}{\mu^{p-1}})^{p^i}b_{i+1}+b_i= 0 \quad\text{ for } i\ge 1\\
&(\frac{p}{\mu^{p-1}})^pb_1+ b_0= a^p     \\
&(\frac{p}{\mu^{p-1}})^pb_0= pa-j\mu.
\end{align*}
 So $b_i= 0$ if $i\ge 1$, $b_0=a^p$ and $pa-j\mu=\frac{p}{\mu^{p-1}}a^p$.
\end{proof}
\begin{ex}\label{ex:ext Z/pZ by Z/p }Let us suppose $\mu=\lb=\lb_{(1)}$.
Then $G_{\lb_{(1)}}\simeq \Z/p\Z$. By \ref{ex:Z/p^2Z SS},
\ref{cor:clE se lb divide mu} and \ref{teo:ext1(glx,gmx)} we have
that $$ \{(k\eta,k)|k\in\Z/p\Z\}\In
\Phi_{\lb_{(1)},\lb_{(1)}}\simeq rad_{p,\lb_{(1)}}(<1+\lb_{(1)}
S>).$$ On the other hand  by \ref{teo:ext1(glx,gmx)} and
\ref{cor:modelli di
Z/p^2Z} it follows that 
$rad_{p,\lb_{(1)}}(<1+\lb_{(1)} S>)\simeq
H^2_0(\Z/p\Z,\Z/p\Z)\simeq \Z/p\Z$. Therefore $
\{(k\eta,k)|k\in\Z/p\Z\}\simeq rad_{p,\lb_{(1)}}(<1+\lb_{(1)}
S>)$. 
\end{ex}

We now concentrate on to the case $v(\mu)\ge v(\lb)$, which is the
unique case, as proved in \ref{prop:rad_p}, where extensions of
$\gmx{1}$ by $\glx{1}$ could be models of $\Z/p^2\Z$, as group
schemes. Our task is to find explicitly all the solutions
$(a,j)\in (R/\lb R)^{\fr}$ of the equation $pa-j\mu=a^p\in R/\lb^p
R$. By \ref{cor:clE se lb divide mu} this means  finding
explicitly all the extensions  of type $\clE^{(\mu,\lb;F,j)}$. Let
us consider the restriction map
$$
 r: \{\clE^{(\mu,\lb;F,j)}
\} \too \Ext^1_K(\mu_p,\mu_p)\simeq \Z/p\Z.
$$
We remark that it coincides with the projection
\begin{align*}
p_2:\bigg\{(a,j)\in (R/\lb R)^{\fr}\times\Z/p\Z& \text{ such that
} pa-j\mu=\frac{p}{\mu^{p-1}}a^p\in R/\lb^p R\bigg\}\too \Z/p\Z
\end{align*}
So there is an extension of $\gmx{1}$ by $\glx{1}$ which is a
model of $(\Z /p^2\Z)_K$ if and only if $p_2$ is surjective. First
of all we  describe explicitly the kernel of the above map.
\begin{lem}\label{lem:ker p2}
We have
\begin{align*}
\ker p_2= \bigg\{(a,0)\in R/\lb R\times&\Z/p\Z \text{ s. t., for
any
  lifting }\tilde{a}\in R,\\ &
pv(\tilde{a})\ge\max\{pv(\lb)+(p-1)v(\mu)-v(p),v(\lb)\}\bigg\}
\end{align*}
In particular $p_2$ is injective if and only if $v(\lb)\le 1$ or
$v(p)-(p-1)v(\mu)<p$.
\end{lem}
\begin{proof}
Let $(a,0)\in \ker p_2\cap R/\lb R\times \Z/p\Z $. By the
definitions we have that
$$
pa=\frac{p}{\mu^{p-1}}a^p\in R/\lb^p R \quad \text{ and } \quad
a^p=0\in R/\lb R.
$$
Let $\tilde{a}\in R$ be a lift of $a$. Since $v(\mu)\ge v(\lb)$,
if $a\neq 0$ then $v(\tilde{a})<v(\mu)$. Hence

\begin{equation}\label{eq:valutazioni diverse}
v(p)+v(\tilde{a})>pv(\tilde{a})+v(p)-(p-1)v(\mu)
\end{equation}

 Therefore
$$
pa=\frac{p}{\mu^{p-1}}a^p\in R/\lb^p R
$$
if and only if
$$
\frac{p}{\mu^{p-1}}a^p=0\in R/\lb^p R,
$$
if and only if
$$
pv(\tilde{a})+v(p)-(p-1)v(\mu)\ge pv(\lb).
$$
We remark that $a^p=0\in R/\lb R$ means $pv(\tilde{a})\ge v(\lb)$.
So we have  proved the first assertion. Now if $v(\lb)\le 1$ or
$v(p)-(p-1)v(\mu)<p$ it is easy to see that there are no nonzero
elements in $\ker p_2$. While if $v(\lb)>1$ and
$v(p)-(p-1)v(\mu)\ge p$,  take $a\in R/\lb R$ with a lifting
$\tilde{a}\in R$ of valuation $v(\lb)-1$. Therefore
$$
p(v(\lb)-1)\ge \max\{pv(\lb)-v(p)+(p-1)v(\mu),v(\lb)\}.
$$
Hence $(a,0)\in \ker p_2$.
\end{proof}

We remark that $\ker p_2$ depends only on the valuations of $\mu$
and $\lb$. So we can easily compute $\Phi_{\mu,\lb}$, too.

\begin{prop}\label{cor:p2 surjective}
Let us suppose $p>2$. Let $\mu, \lb \in R\setminus \{0\}$ be with
$v(\lb_{(1)})\ge v(\mu)\ge v(\lb)$.
\begin{itemize}
\item[a)] If $v(\mu)<pv(\lb)$ then $p_2$ is surjective if and only
if $pv(\mu)-v(\lb)\ge v(p)$. And, if $p_2$ is surjective, $
 \Phi_{\mu,\lb}
$ is isomorphic to the group
$$
\{(j \eta \frac{\mu}{\lb_{(1)}}+ \alpha,j)| (\alpha,0)\in
\ker(p_2) \text{ and } j\in \Z/p\Z\}
$$
For the definition of $\eta$ see \ref{ex:Z/p^2Z SS}. 
\item[b)] If $v(\mu)\ge pv(\lb)$ then $p_2$ is surjective and $
\Phi_{\mu,\lb}$ is isomorphic to
$$
\{(\alpha,j)| (\alpha,0)\in \ker(p_2) \text{ and } j\in
\Z/p\Z\}\simeq \ker p_2\times \Z/p \Z.
$$
\item[c)] If $p_2$ is not surjective then $p_2$ is the zero
morphism. So $
 \Phi_{\mu,\lb}=\ker p_2
$.
\end{itemize}
\end{prop}
\begin{rem}\label{rem:valutazione di elementi di Phi1}Let us suppose $v(\mu)<pv(\lb)$. Let $(b,j)\in \Phi_{\mu,\lb}$ with $j\neq 0$. By \ref{rem:0 in Phi}, then $b\neq 0$. Let $\tilde{b}\in
R$ be any of its lifting. Then
$v(\tilde{b})=v(\eta\frac{\mu}{\lb_{(1)}})=v(\mu)-\frac{v(p)}{p}$.
Indeed, by the theorem, we have
$\tilde{b}=\eta\frac{\mu}{\lb_{(1)}}+\alpha$ for some $\alpha \in
R/\lb R$ with
$v(\tilde{\alpha})>v(\eta\frac{\mu}{\lb_{(1)}})=v(\mu)-\frac{v(p)}{p}$,
where $\tilde{\alpha}\in R$ is any lifting of $\alpha$.
\end{rem}
\begin{proof}
\begin{itemize}
\item[a)] First, we suppose that $p_2$ is surjective. This is
equivalent to saying that
\begin{equation}\label{eq:soluzione pa-jmu...}
pa-j\mu=\frac{p}{\mu^{p-1}}a^p\in R/\lb R
\end{equation}
has a solution $a\in (R/\lb R)^{\fr}$ if $j\neq 0$. Since
$v(\mu)<v(p)$, by \eqref{eq:valutazioni diverse} it follows that
\begin{equation}\label{eq:soluzione pa=... segue che...}
v(\mu)=v(p)-(p-1)v(\mu)+pv(\tilde{a}),
\end{equation}
with $\tilde{a}\in R$ a lifting of $a$. Since $a\in (R/\lb
R)^{\fr}$ we have $pv(\tilde{a})\ge v(\lb)$. Hence, by
\eqref{eq:soluzione pa=... segue che...}, $pv(\mu)- v(\lb)\ge
v(p)$.

Conversely let us suppose that $pv(\mu)-v(\lb)\ge v(p)$. We know
by \ref{ex:Z/p^2Z SS} and \ref{cor:modelli di Z/p^2Z} that
$$
p\eta-\lb_{(1)}=\frac{p}{\lb_{(1)}^{p-1}}\eta^p\in R/ \lb_{(1)}^p
R.
$$
We recall that $v(\eta)=v(\lb_{(2)})$.  Since $p
v(\lb_{(1)})-v(\lb_{(1)})+\mu\ge p\mu\ge p\lb$, if we divide the
above equation by $\frac{\lb_{(1)}}{\mu}$ we obtain
\begin{equation*}
p\eta \frac{\mu}{\lb_{(1)}}-\mu=
\frac{p}{\mu^{(p-1)}}(\frac{\mu}{\lb_{(1)}}\eta)^p\in R/\lb^p R.
\end{equation*}
We remark that $\eta \frac{\mu}{\lb_{(1)}}\in (R/\lb R)^{\fr}$,
since, by hypothesis, $v(\big(\eta
\frac{\mu}{\lb_{(1)}}\big)^p)=pv(\mu)-v(p)\ge v(\lb)$. Clearly
$j\eta \frac{\mu}{\lb_{(1)}}$ is a solution of \eqref{eq:soluzione
pa-jmu...} for any $j\in \Z/p\Z$.

In particular it follows that, if $p_2$ is surjective,
$\Phi_{\mu,\lb}$ is isomorphic to the group
$$
\{(j \eta \frac{\mu}{\lb_{(1)}}+ \alpha,j)| (\alpha,0)\in
\ker(p_2) \text{ and } j\in\ \Z/p\Z\}
$$
\item[b)] If $v(\mu)\ge pv(\lb)$ then we have that $\mu=0\in
R/\lb^p R$.  We remark that $(0,j)\in \Phi_{\mu,\lb}$. This
implies that  $p_2$ is surjective and that  $(\alpha,j)\in R/\lb
R\times \Z/p\Z \cap \Phi_{\mu,\lb}$ if and only if $(\alpha,0)\in
\ker(p_2)$.

\item[c)] Since $p_2$ is a morphism of groups with target $\Z/p\Z$
then the image of $p_2$ is a subgroup of $\Z/p\Z$. Then the image
of $p_2$  is trivial or it is equal to $\Z/p\Z$. The assertion
follows.
\end{itemize}
\end{proof}
\begin{ex}\label{ex:Ext 1 con mu=lb_(1)}
Let us suppose 
$v(\mu)=v(\lb_{(1)})$, i.e. $G_{\mu,1}\simeq \Z/p\Z$. For
simplicity we will suppose $\mu=\lb_{(1)}$. Then $p_2$ is an
isomorphism. Indeed in this case $\ker(p_2)=0$ by \ref{lem:ker p2}
and it is surjective by \ref{cor:p2 surjective}(a)-(b). This means
that, in this case, any extension $\clE^{(\lb_{(1)},\lb;F,j)}$ is
uniquely determined by the induced extension  over $K$. Let us now
consider the map
$$
\Ext^1(G_{\lb_{(1)},1},G_{\lb_{(1)},1})\too
\Ext^1(G_{\lb_{(1)},1},G_{\lb,1})
$$
induced by the map $\Z/p\Z\simeq G_{\lb_{(1)},1}\too G_{\lb,1}$
given by $S\mapsto \frac{\lb_{(1)}}{\lb}S$. It is easy to see that
$\clE^{(\lb_{(1)},\lb;F,j)}$ is the image of
$\clE^{(\lb_{(1)},\lb_{(1)};E_p(\eta S),j)}$ through the above
map. Indeed from the above proposition we have that $F(S)\equiv
E_p(\eta S)\mod \lb$. We remark that if $pv(\lb)\le v(\lb_{(1)})$
then $\eta\equiv 0\mod \lb$, indeed in such a case $v(\lb)\le
v(\lb_{(2)})=v(\eta)$.
\end{ex}

\subsection{Classification of models of $(\Z/p^2\Z)_K$}
By the previous paragraphs we have a classification of extensions
of $G_{\mu,1}$ by $G_{\lb,1}$ whose generic fibre is isomorphic,
as group scheme, to $\Z/p^2\Z$. But this classification is too
fine for our tasks. We want here to forget the structure of
extension. We are only interested in the group scheme structure.
We observe that it can happen that two non isomorphic extensions
are isomorphic as group schemes. We here study when it happens.

First of all we recall what the model maps between models of $\Z/p
\Z$  are.
 Let us suppose  $\rho,\tilde{\rho}\in R$ with $v(\rho),v(\tilde{\rho})\le v(\lb_{(1)})$. Since $G_{\rho,1}$ is flat over $R$,
by \ref{lem:morfismo gen nullo=nullo} it follows that the
restriction map
$$
\Hom_{gr}(G_{\rho,1},G_{\tilde{\rho},1})\too
\Hom_{gr}(({G_{\rho,1}})_K,{(G_{\tilde{\rho},1})}_K)\simeq \Z/p\Z
$$
in an injection. It  follows easily by \eqref{eq:Hom(gmu1,gl)}
that
$$
\Hom_{gr}(G_{\rho,1},G_{\tilde{\rho},1})=
\left\{
\begin{array}{ll}
    \Z/p\Z, & \hbox{ if $v(\rho)\ge v(\tilde{\rho})$;} \\
0, & \hbox{if $v(\rho)<v(\tilde{\rho}$)},
\\
\end{array}%
\right.
$$
where, in the first case, the morphisms are given by $S\longmapsto
\frac{(1+\rho S)^r-1}{\tilde{\rho}}$ with $r\in \Z/p\Z$. We remark
that, if $v(\rho)=v(\tilde{\rho})$ and  $r\neq 0$, these morphisms
are isomorphisms. 

We now recall that by \ref{lem:modelli di Z/p^2 Z sono
estensioni}, \ref{cor:modelli di Z/p^2Z}, \ref{cor:clE se lb
divide mu} and \ref{prop:rad_p} any model of $(\Z/p^2\Z)_K$ is of
the form $\clE^{(\mu,\lb;F,j)}$ such that $j\neq 0$,
$v(\lb_{(1)})\ge v(\mu)\ge v(\lb)$ and
$F(S)=\sum_{i=0}^{p-1}\frac{a^i}{i!}S^i$ with $(a,j)\in
\Phi_{\mu,\lb}$. See \ref{cor:p2 surjective} for the explicit
description of $\Phi_{\mu,\lb}$.
For $i=1,2$ let us consider $\clE^{(\mu_i,\lb_i;F_i,j_i)}$, models
of $(\Z/p^2\Z)_K$. First of all we remark that there is an
injection
$$
r_K:\Hom(\clE^{(\mu_1,\lb_{1},F_1,j_1)},
\clE^{(\mu_2,\lb_{2},F_2,j_2)})\too
\Hom_K(\clE_{j_1,K},\clE_{j_2,K})
$$
given by
$$
f\longmapsto {(\alpha_{\mu_2,{\lb_2}})}_K\circ f_K \circ
{(\alpha_{\mu_1,\lb_1})}^{-1}_K
$$
See  \eqref{eq:morfismo Z/p2Z to E} for the definition
$\alpha_{\mu,\lb}$.  We recall that
$$
\Hom({\clE_{j_1}},{\clE_{j_2}})\simeq
\Hom_K(\clE_{j_1,K},\clE_{j_2,K}).
$$
and the  elements are the morphisms
$$
\psi_{r,s}:{\clE_{j_1}}\too \clE_{j_2},
$$
which, on the level of Hopf algebras, are given by
\begin{align}\label{eq:morfismi tra cle j}
S_1&\longmapsto S_1^{\frac{rj_1}{j_2}}\\
S_2&\longmapsto S_1^sS_2^r,
\end{align}
for some $r\in \Z/p\Z$ and $s\in \Z/p\Z$.  Moreover the map
\begin{align*}
\Hom({\clE_{j_1}},{\clE_{j_2}})&\too \Z/p^2\Z\\
\psi_{r,s}&\longmapsto r+\frac{p}{j_1}s
\end{align*}
is an isomorphism. So $\Hom(\clE^{(\mu_1,\lb_{1},F_1,j_1)},
\clE^{(\mu_2,\lb_{2},F_2,j_2)})$ is a subgroup of $\Z/p^2\Z$
through the map $r_K$. We remark that the unique nontrivial
subgroup of $\Hom({\clE_{j_1}},{\clE_{j_2}})$ is
$\{\psi_{0,s}|s\in \Z/p\Z\}$. Finally we have that any morphism
$\clE^{(\mu_1,\lb_{1},F_1,j_1)}\too
\clE^{(\mu_2,\lb_{2},F_2,j_2)})$ is given by
\begin{equation}\label{eq:def f}
\begin{aligned}
S_1&\too \frac{(1+\mu_1 S_1)^{\frac{rj_1}{j_2}}-1}{{\mu_2}} \\
S_2&\too \frac{(F_1(S_1)+\lb_{1} S_2)^{r}(1+\mu_1
S_1)^s-F_2(\frac{(1+\mu_1
S_1)^{\frac{rj_1}{j_2}}-1}{{\mu_2}})}{{\lb_{2}}},
\end{aligned}
\end{equation}
for some $r,s\in \Z/p\Z$. With abuse of notation we call it
$\psi_{r,s}$. We remark that the morphisms
$\psi_{r,s}:\clE^{(\mu_1,\lb_1,F_1,j_1)}\too
\clE^{(\mu_2,\lb_2,F_2,j_2)}$ which are model maps correspond, by
\eqref{eq:morfismi tra cle j}, to $r\neq 0$. In such a case
$\psi_{r,s}$ is a morphism of extensions, i.e. there exist
morphisms $\psi_1:G_{\lb,1}\too G_{\lb,2}$ and
$\psi_2:G_{\mu,1}\too G_{\mu,2}$ such that
\begin{equation}\label{eq:model maps are morfismi di estensioni}
\begin{aligned}
\xymatrix@1{0\ar[r]&{G_{\lb_1,1}}\ar[r]^{}\ar[d]^{\psi_1}&{\clE^{(\mu_1,\lb_1,F_1,j_1)}}\ar^{\psi_{r,s}}[d]\ar[r]&{G_{\mu_1,1}}\ar[r]\ar[d]^{\psi_2}&0\\
                     0\ar[r]&{G_{\lb_2,1}}\ar[r]^{}&\clE^{(\mu_2,\lb_2,F_2,j_2)}\ar[r]^{}&{G_{\mu_2,1}}\ar[r]&0}
\end{aligned}
\end{equation}
commutes. More precisely $\psi_1$ is given by $S\mapsto
\frac{(1+\lb_1 S)^{r}-1}{{\lb_2}}$ and $\psi_2$ by $S\mapsto
\frac{(1+\mu S_1)^{\frac{rj_1}{j_2}}-1}{{\mu_2}}$. 

We now calculate $\Hom(\clE^{(\mu_1,\lb_1,F_1,j_1)},
\clE^{(\mu_2,\lb_{2},F_2,j_2)})$.
\begin{prop}\label{lem:abbasso valutazione con blow-up} 
 For $i=1,2$, if $F_i(S)=E_p(a_i S)=\sum_{k=0}^{p-1}\frac{a_i^k}{k!}S^i$ and
${\scE}_i=\clE^{(\mu_i,\lb_{i};F_i,j_i)}$ are models of
$(\Z/p^2\Z)_K$  we have
$$
\Hom(\scE_1,\scE_2)=
\left\{%
\begin{array}{ll}
    0, & \hbox{\text{if $v(\mu_1)<v(\lb_2)$;}} \\
    \{\psi_{r,s}\}\simeq \Z/p^2\Z, & \hbox{\text{if $v(\mu_2)\le v(\mu_1)$, $v(\lb_2)\le v(\lb_1)$}}\\
    &  \text{and $a_1\equiv \frac{j_1}{j_2}\frac{\mu_1}{\mu_2}a_2\mod{\lb_2}$;}\\
\{\psi_{0,s}\}\simeq \Z/p\Z, & \hbox{\text{otherwise}}.
\end{array}%
\right.
$$

\end{prop}
\begin{proof}
It is immediate to see that $\psi_{0,s}\in \Hom(\scE_1,\scE_2)$,
with $s\neq 0$,  if and only if $v(\mu_1)\ge v(\lb_2)$. We now see
conditions for the existence of $\psi_{r,s}$ with $r\neq 0$.
 If it exists, in particular, we have two  morphisms
 $G_{\mu_1,1}\too G_{\mu_2,1}$ and $G_{\lb_1,1}\too G_{\lb_2,1}$.
This implies $v(\mu_1)\ge v(\mu_2)$ and $v(\lb_1)\ge v(\lb_2)$.
Moreover we have that
$$
F_1(S_1)^r(1+\mu_1 S_1)^s=F_2(\frac{(1+\mu_1
S_1)^{\frac{rj_1}{j_2}}-1}{{\mu_2}})\in
\Hom({G_{\mu_1,1}}_{|S_{\lb_{2}}},{\gm}_{|S_{\lb_{2}}}).
$$
Since $v(\mu_1)\ge v(\mu_2)\ge v(\lb_2)$,  we  have
\begin{equation}\label{eq:F_1^r=F_2^s...}
F_1(S_1)^r=F_2(\frac{(1+\mu_1
S_1)^{\frac{rj_1}{j_2}}-1}{{\mu_2}})\in
\Hom({G_{\mu_1,1}}_{|S_{\lb_{2}}},{\gm}_{|S_{\lb_{2}}}).
\end{equation}
If we define the morphism of groups
\begin{align*}
[\frac{\mu_1}{\mu_2}]^*:\Hom({G_{\mu_2,1}}_{|S_{\lb_{2}}},{\gm}_{|S_{\lb_{2}}})&\too
\Hom({G_{\mu_1,1}}_{|S_{\lb_{2}}},{\gm}_{|S_{\lb_{2}}})\\
F(S_1)&\longmapsto F(\frac{\mu_1}{\mu_2}S_1)
\end{align*}
then
\begin{align*}
F_2(\frac{(1+\mu_1
S_1)^{\frac{rj_1}{j_2}}-1}{\mu_2})&=[\frac{\mu_1}{\mu_2}]^*\bigg(F_2\bigg(\frac{(1+\mu_1
S_1)^{\frac{rj_1}{j_2}}-1}{\mu_1}\bigg)\bigg)\\&=[\frac{\mu_1}{\mu_2}]^*(F_2(S_1))^{\frac{rj_1}{j_2}}\\
&=F_2(\frac{\mu_1}{\mu_2}(S_1))^{\frac{rj_1}{j_2}}.
\end{align*}
Therefore we have
\begin{equation}\label{eq:F1=F=2}
F_1(S_1)^r=(F_2(\frac{\mu_1}{\mu_2}S_1))^{\frac{rj_1}{j_2}}\in
\Hom({G_{\mu_1,1}}_{|S_{\lb_{2}}},{\gm}_{|S_{\lb_{2}}}).
\end{equation}

 Every element
of $\Hom({G_{\mu_1,1}}_{|S_\lb},{\gm}_{|S_\lb}) $ has order $p$.
Let $t$ be an inverse for $r$  modulo $p$. Then raising the
equality to the $t^{th}$-power we obtain
$$
F_1(S_1)=(F_2(\frac{\mu_1}{\mu_2}S_1))^{\frac{j_1}{j_2}}\in\Hom({G_{\mu_1,1}}_{|S_{\lb_2}},{\gm}_{|S_{\lb_2}}).
$$
By \ref{cor: hom gmu gm se lb divide mu} this means
$$
a_1\equiv {\frac{j_1}{j_2}}\frac{\mu_1}{\mu_2}a_2\mod \lb_2.
$$

 It is conversely clear
that,  if  $v(\mu_1)\ge v(\mu_2)$, $v(\lb_{1})\ge v({\lb_{2}})$
and
$$
F_1(S_1)=(F_2(\frac{\mu_1}{\mu_2}S_1))^{\frac{j_1}{j_2}}\in
\Hom({G_{\mu_1,1}}_{|S_{\lb_{2}}},{\gm}_{|S_{\lb_{2}}}),
$$
then \eqref{eq:def f} defines a morphism of group schemes. 

 \end{proof}
We have the following result which gives a criterion to determine
the class of isomorphism, as a group scheme, of an extension of
type $\clE^{(\mu,\lb;F,j)}.$
\begin{cor}\label{cor:iso tra modelli}
For $i=1,2$, let $F_i(S)=E_p(a
S)=\sum_{k=0}^{p-1}\frac{a_i^k}{k!}S^k$ and let
$\scE_i=\clE^{(\mu_i,\lb_{i};F_i,j_i)}$ be models of
$(\Z/p^2\Z)_K$. Then they are isomorphic if and only if
$v(\mu_1)=v(\mu_2)$, $v(\lb_1)=v(\lb_2)$ and $a_1\equiv
    \frac{j_1}{j_2}\frac{\mu_1}{\mu_2}a_2\mod{\lb_2}$. Moreover if
    it happens then any model map between them is an isomorphism.
\end{cor}
\begin{proof}
By the proposition we have that  a model map
$\psi_{r,s}:\clE^{(\mu_1,\lb_1,F_1,j_1)}\too
\clE^{(\mu_2,\lb_2,F_2,j_2)}$  exists if and only if $v(\mu_1)\ge
v(\mu_2)$, $v(\lb_1)\ge v(\lb_2)$ and $a_1\equiv
\frac{j_1}{j_2}\frac{\mu_1}{\mu_2}a_2\mod{\lb_2}$. It is a
morphism of extensions as remarked before the proposition.
    Let us consider the commutative diagram \eqref{eq:morfismi tra cle j}. Then $\psi_{r,s}$ is an isomorphism if and only if $\psi_i$ is an isomorphism for $i=1,2$.
    By the discussion made at the beginning of this section  this is equivalent to requiring $v(\mu_1)=v(\mu_2)$ and $v(\lb_{1})=v(\lb_{2})$.
    This also proves the last
 assertion. 
\end{proof}
\vspace{.7cm} We  remarked that if $v(\mu_1)=v(\mu_2)$ and
$v(\lb_1)=v(\lb_2)$ then
$$
\Ext^1(G_{\mu_1,1},G_{\lb_1,1})\simeq\Ext^1(G_{\mu_2,1},G_{\lb_2,1}).
$$
The following is a more precise statement for extensions of type
$\clE^{(\mu,\lb;F,j)}$.
\begin{cor}\label{cor:scelgo mu'} Let $\clE^{(\mu_1,\lb_1;E_p(a S),j)}\in \ext^{1}(G_{\mu_1,1},G_{\lb_{1},1})$ be a model of $\Z/p^2\Z$. 
Then for any $\mu_2,\lb_2$ such that $v(\mu_1)=v(\mu_2)$ and $v(\lb_1)=v(\lb_2)$ we 
have
$$
\clE^{(\mu_1,\lb_1;E_p(aS
),j)}\simeq\clE^{({\mu_2},\lb_2;E_p(\frac{a}{j}\frac{\mu_2}{\mu_1}S),1)}
$$
as group schemes.
\end{cor}
\begin{proof} Firstly we  prove that there exists the group scheme $\clE^{({\mu_2},\lb_2;E_p(\frac{a}{j}\frac{\mu_2}{\mu_1}S),1)}$.
By \ref{cor:clE se lb divide mu} we have 
that $a\in (R/\lb R)^{\fr}$ and
\begin{equation}\label{eq:equazione per bbf a}
pa-j\mu_1=\frac{p}{\mu_1^{p-1}}a^p\mod\lb_1^p.
\end{equation}
Then, multiplying \eqref{eq:equazione per bbf a} by
$\frac{\mu_2}{\mu_1}\frac{1}{j}$, we  have
$$
p\frac{a \mu_2}{j\mu_1}-\mu_2\equiv
\frac{p}{\mu_2^{p-1}}\bigg(\frac{a{\mu_2}}{ j\mu_1}\bigg)^p\mod
{\lb_2}^p.
$$
Hence
$\clE^{({\mu_2},\lb_2;E_p(\frac{a}{j}\frac{\mu_2}{\mu_1}S),1)}$ is
a group scheme (see again  \ref{cor:clE se lb divide mu}). 
Then by the above proposition we can conclude that
$$
\clE^{(\mu_1,\lb_1;E_p(a S),j
)}\simeq\clE^{({\mu_2},\lb_2;E_p(\frac{a}{j}\frac{\mu_2}{\mu_1}S),1)}
$$
as group schemes.

\end{proof}
\begin{ex}\label{ex:equazione per eta_pi}
    Let $\mu,\lb\in R$ be such that $v(\mu)=v(\lb)=v(\lb_{(1)})$. We now
    want to describe $\Z/p^2\Z$ as
    $\clE^{(\mu,\lb;F,1)}$. We recall that we defined $$
\eta=\sum_{k=1}^{p-1}\frac{(-1)^{k-1}}{k}\lb_{(2)}^k
$$  By \ref{ex:Z/p^2Z SS} and the previous corollary we have
$$
\Z/p^2\Z\simeq \clE^{(\mu,\lb;E_p(\eta\frac{\mu}{\lb_{(1)}}S),1)}.
$$
    \end{ex}
\vspace{.3cm}

We conclude the section with the complete classification of
$(\Z/p^2\Z)_K$-models. The following theorem summarizes the above
results.

 \begin{thm}\label{cor:modelli di Z/p^2 Z sono cosi'2}
Let us suppose $p>2$. Let $G$ be a finite and flat $R$-group
scheme such that $G_K\simeq (\Z/p^2\Z)_K$. Then $G\simeq
\clE^{(\pi^{m},\pi^{n};E_p(a S),1)}$ for some $v(\lb_{(1)})\ge
m\ge n\ge 0$ and $(a,1)\in \Phi_{\pi^m,\pi^n}$. Moreover $m,n$ and
$a\in R/\pi^{n} R$ are unique.
\end{thm}
\begin{rem} The explicit description of the set $\Phi_{\pi^m,\pi^n}$ has been given in \ref{cor:p2 surjective} and \ref{lem:ker p2}.
\end{rem}
 \begin{proof} By \ref{lem:modelli di Z/p^2 Z sono estensioni}, \ref{cor:modelli
di Z/p^2Z}, \ref{prop:rad_p}, \ref{cor:clE se lb divide mu} and
\ref{cor:scelgo mu'}
 any model of $(\Z/p^2\Z)_K$ is of type $\clE^{(\pi^m,\pi^n;E_p(a S),1)}$ with $m\ge n$ and $(a,1)\in\Phi_{\pi^m,\pi^n}$. By \ref{cor:iso tra modelli}, it follows that,
$$
 \clE^{(\pi^{m_1},\pi^{n_1},E_p(a_1 S),1)}\simeq
\clE^{({\pi^{m_2}},{\pi^{n_2}},E_p(a_2 S),1)}
$$
as group schemes  if and only if $m_1=m_2$, $n_1=n_2$ and $a_1=a_2\in R/\pi^{n_1} R$.  
\end{proof}

\section{Reduction on the special fiber of the models of $(\Z/p^2\Z)_K$}
In the following we study the special fibers of the extensions of
type $\clE^{(\lb,\mu,;F,j)}$ with $v(\mu)\ge v(\lb)$. In
particular, by \ref{prop:rad_p}, this includes the  extensions
which are
models of $(\Z/p^2\Z)_K$ as group schemes. 
We study separately the different cases which can occur.


%
%
\subsection{Case $\mathbf{v(\mu)=v(\lb)=0}$}

We have ${(G_{\lb,1})}_k\simeq {(G_{\mu,1})}_k\simeq \mu_p$. The
extensions of type $\clE^{(\lb,\mu,;F,j)}$ are the extensions $
\clE_i$ with $i\in \Z/p\Z$. The
special fibers of the extensions  
$\clE_i$ with $i\in \Z/p\Z$ are clearly $\clE_{i,k}$. 
See also \ref{eq:ext1(mup,mup)}.
\subsection{Case $\mathbf{v(\lb_{(1)})\ge v(\mu)>v(\lb)=0}$}

In such a case we have ${(G_{\lb,1})}_k\simeq \mu_p$.
 It is immediate by the
definitions that any  extension $\clE^{(\mu,\lb;1,j)}$ is trivial
on the special fiber.
\subsection{Case $\mathbf{v(\lb_{(1)})>v(\mu)\ge v(\lb)>0}$}
 Then ${(G_{\mu,1})}_k\simeq (G_{\lb,1})_k\simeq
{\ap}_{,k}$.

 First, we recall some results about extensions of group schemes
of order $p$ over a field $k$. See \cite[III \S 6 7.7.]{DG} for a
reference.

\begin{thm}\label{teo:ext(ap,ap)} Let us suppose that $\ap$ acts trivially on $\ap$ over $k$.
The exact sequence $$ 0\too \alpha_p\to \Ga\on{\fr}\to \Ga\to 0$$
induces the following split exact sequence
$$
0\too \Hom_k (\alpha_p,\Ga)\too \Ext^1(\ap,\ap)\too
\Ext^1(\ap,\Ga)\too 0.
$$
\end{thm}
It is also known that
 $$
\Ext^1(\Ga,\Ga)\simeq H_0^2(\Ga,\Ga)\too H_0^2(\ap,\Ga)\simeq
\Ext^1(\ap,\Ga).
 $$
 is surjective. Since $\Ext^1(\Ga,\Ga)\simeq H_0^2(\Ga,\Ga)$ is
 freely generated as a right $k[\fr]$-module by
 $C_i=\frac{X^{p^i}+X^{p^i}-(X+Y)^{p^i}}{p^i}$ and $D_i=X Y^{p^i}$ for all $i\in \N\setminus\{0\}$, it
 follows that
 $H_0^2(\ap,\Ga)\simeq \Ext^1(\ap,\Ga)$ is freely generated as
right $k$-module by the class of the  cocycle
$C_1=\frac{X^p+X^p-(X+Y)^p}{p}$. So $\Ext^1(\ap,\Ga)\simeq k$.

 Moreover it is easy to see that $\Hom_k (\alpha_p,\Ga)\simeq k$.
 The morphisms are given by $T\mapsto a T$ with $a\in k$.
By these remarks we have that the isomorphism
$$
\Hom_k (\alpha_p,\Ga)\times \Ext^1(\ap,\Ga)\too \Ext^1(\ap,\ap),
$$
deduced from \ref{teo:ext(ap,ap)}, is given by
$$
(\beta,\gamma C_1)\mapsto E_{\beta,\gamma}.
$$
The extension $E_{\beta,\gamma}$ is so defined:
$$
E_{\beta,\gamma}=\Sp(k[S_1,S_2]/(S_1^p,S_2^p-\beta S_1))
$$

\begin{enumerate}
    \item law of multiplication
    \begin{align*}
    S_1\longmapsto &S_1\pt 1+1\pt S_1\\
    S_2\longmapsto &S_2\pt 1 + 1\pt S_2 + \gamma \frac{S_1^p\pt 1+ 1\pt S_1^p-(S_1\pt 1+1\pt S_1)^p}{p}
    \end{align*}
    \item unit
    \begin{align*}
    &S_1\longmapsto 0\\
    &S_2\longmapsto 0
    \end{align*}
    \item inverse
    \begin{align*}
    &S_1\longmapsto -S_1\\
    &S_2\longmapsto -S_2
\end{align*}
\end{enumerate}
It is clear that all such extensions are commutative. In
\cite[4.3.1]{SS4} the following result  was proved.
\begin{prop}Let $\lb,\mu\in \pi R\setminus\{0\}$. Then $[\clE^{(\mu,\lb;E_{p}(\bbf{a},\mu,S))}_k]\in H_0^2({\Ga}_{,k},{\Ga}_{,k})$
coincides with the class of
$$
\sum_{k=1}^{\infty}\frac{(\fr-{[\mu^{p-1}])}(\tilde{\bbf{a}})}{\lb}C_k,
$$
where $\tilde{\bbf{a}}\in \widehat{W}(R)$ is a lifting of
$\bbf{a}\in \widehat{W}(R/\lb R)$.
\end{prop}
We deduce the following corollary about the extensions of
$\alpha_p$ by $\Ga$.
\begin{cor}\label{cor:ext(ap,ga)}Let $\lb,\mu\in \pi R\setminus\{0\}$. Then $[\tilde{\clE}^{(\mu,\lb;E_{p}(\bbf{a},\mu,S))}_k]\in H_0^2({\ap}_{,k},{\Ga}_{,k})$
coincides with the class of
$$
\frac{(\fr-{[\mu^{p-1}])}(\tilde{\bbf{a}})}{\lb}C_1,
$$
where $\tilde{\bbf{a}}  \in \widehat{W}(R)$ is a lifting of
$\bbf{a}\in \widehat{W}(R/\lb R)$.
\end{cor}
\begin{proof}
This follows from the fact that
${\clE}^{(\mu,\lb;E_{p}(\bbf{a},\mu,S))}_k\mapsto
\tilde{\clE}^{(\mu,\lb;E_{p}(\bbf{a},\mu,S))}_k$ through the map
$$ \Ext^1(\Ga,\Ga)\simeq H_0^2(\Ga,\Ga)\too H_0^2(\ap,\Ga)\simeq
\Ext^1(\ap,\Ga).
 $$
\end{proof}

Let us take an extension $\clE^{(\mu,\lb;E_p(a S),j)}$. Let
$\tilde{a}$ be a lifting of $a$. We have that on the special fiber
this extension is given as a scheme by
$$
\clE^{(\mu,\lb;E_p(a
S),j)}_k=\Sp(k[S_1,S_2]/(S_1^p,S_2^p-(-\frac{(\sum_{i=0}^{p-1}\frac{\tilde{a}^i}{i!}S^i)^p(1+\mu
S_1)^{-j}-1}{\lb^p}))).
$$

By \ref{cor:clE se lb divide mu} we know that
$$
pa-j\mu-\frac{p}{\mu^{p-1}}a^p=0\in R/\lb^p R
$$
In the proof of the same corollary we have seen that
$$
p[a]-j[\mu]-[\frac{p}{\mu^{p-1}}a^p]-V([a^p])=[pa-j\mu-\frac{p}{\mu^{p-1}}a^p]\in
\widehat{W}(R/\lb^p R).
$$
 By the definitions we
have the following equality in
$\Hom(\gmu_{|S_{\lb^p}},{\gm}_{|S_{\lb^p}})$
$$
\xi^0_{R/\lb^p
R}(p[a]-j[\mu]-[\frac{p}{\mu^{p-1}}a^p]-V([a^p]))=E_p(aS_1)^p(1+\mu
S_1)^{-j}\bigg(E_p\bigg(a^p \bigg(\frac{(1+\mu S_1)^p-1}{\mu }
\bigg)\bigg)\bigg)^{-1}.
$$
Moreover we have
$$ \xi^0_{R/\lb^p
R}([pa-j\mu-\frac{p}{\mu^{p-1}}a^p])=E_p((pa-j\mu-\frac{p}{\mu^{p-1}}a^p)
S_1)$$

So we have that
\begin{align*}
(\sum_{i=0}^{p-1}\frac{\tilde{a}^i}{i!}S_1^i)^p(1+\mu
 S_1)^{-j}-1&\equiv\sum_{i=0}^{p-1}\frac{(pa-j\mu-\frac{p}{\mu^{p-1}}a^p)^iS_1^i}{i!}-1\\&\equiv 0 \mod
 \lb^p \bigg(R[S_1]/(\frac{(1+\lb S_1)^p-1}{\lb^p})\bigg).
\end{align*}
Hence
\begin{align*}
\frac{(\sum_{i=0}^{p-1}\frac{\tilde{a}^i}{i!}S_1^i)^p(1+\mu
 S_1)^{-j}-1}{\lb^p}&\equiv\frac{\sum_{i=0}^{p-1}\frac{(pa-j\mu-\frac{p}{\mu^{p-1}}a^p)^i}{i!}S_1^i-1}{\lb^p}\\
 &\equiv \frac{(pa-j\mu-\frac{p}{\mu^{p-1}}a^p)}{\lb^p}S_1\mod
 \pi.
\end{align*}
On the other hand $\clE^{(\mu,\lb;E_{p}(aS),j)}_k\mapsto
\tilde{\clE}^{(\mu,\lb;E_{p}(aS))}_k$ through the map
$\ext^{1}(\ap,\ap)\to \ext^{1}(\ap,\Ga)$.

Therefore $\clE^{(\mu,\lb;E_{p}(aS),j)}_k\simeq E_{\beta,\gamma}$
with
$\beta=(-\frac{p\tilde{{a}}-j\mu-\frac{p}{\mu^{p-1}}\tilde{{a}}^p}{\lb^p}
\mod \pi)$ and $\gamma=(\frac{\tilde{{a}}^p}{\lb}\mod \pi)$.
We have so proved 
the following result.
\begin{prop}
Let $\lb,\mu\in \pi R$ be such that $v(\lb)\le
v(\mu)<v(\lb_{(1)})$. Then $[\clE^{(\mu,\lb;E_{p}(aS),j)}_k]\in
\Ext_k^1({\ap},{\ap})$ coincides with the class of
$$
\bigg(-\frac{p\tilde{{a}}-j\mu-\frac{p}{\mu^{p-1}}\tilde{{a}}^p}{\lb^p},\frac{\tilde{{a}}^p}{\lb}
C_1\bigg),
$$
where $\tilde{{a}}\in R$ is a lifting of $a\in R/\lb R$.
\end{prop}
\subsection{Case $\mathbf{v(\lb_{(1)})=v(\mu)> v(\lb)>0}$}

In this situation we have $${(\gmx{1})}_k\simeq \Z/p\Z \quad
\text{ and } \quad {(\glx{1})}_k\simeq \ap.$$

%
\begin{prop}
Let $\lb,\mu\in \pi R$ be such that $ v(\mu)=v(\lb_{(1)})>v(\lb)$.
Then $\clE^{(\mu,\lb;E_{p}(aS),j)}_k$ is the trivial extension.
\end{prop}
\begin{proof}
We can suppose $\mu=\lb_{(1)}$. From  \ref{ex:Ext 1 con mu=lb_(1)}
it follows that $\clE^{(\lb_{(1)},\lb;F,j)}$ is in the image of
the morphism
$$
\Ext^1(G_{\lb_{(1)},1},G_{\lb_{(1)},1})\too
\Ext^1(G_{\lb_{(1)},1},G_{\lb,1})
$$
induced by the map $\Z/p\Z\simeq G_{\lb_{(1)},1}\too G_{\lb,1}$
given by $S\mapsto \frac{\lb_{(1)}}{\lb}S$.  But this morphism is
the zero morphism on the special fiber. So we are done.
\end{proof}
\subsection{Case $\mathbf{v(\lb_{(1)})=v(\mu)= v(\lb)}$}

 We have
$${(\gmx{1})}_k\simeq \Z/p\Z \quad \text{ and }
{(\glx{1})}_k\simeq \Z/p\Z.$$ For simplicity we will suppose
$\mu=\lb=\lb_{(1)}$. We recall the following result.
\begin{prop}\label{prop:ext(Fp,Fp)}Let suppose that $\Z/p\Z$ acts trivially on $\Z/p\Z$ over $k$. The exact sequence $0\to \Z/p\Z \to \Ga\on{\fr-1}{\to} \Ga\to 0$ induces  the following exact sequence
\begin{align*}
 \Hom_{gr}(\Z/p\Z,\Ga)\simeq k\on{\fr-1}{\too}\Hom_{gr}(\Z/p\Z,\Ga)\simeq
k&\too \Ext_k^1(\Z/p\Z,\Z/p\Z)\too \\
\too\Ext_k^1(\Z/p\Z,\Ga)&\simeq
k\on{\fr-1}{\too}\Ext_k^1(\Z/p\Z,\Ga)\simeq k
\end{align*}
\end{prop}
\begin{proof}
\cite{DG}
\end{proof}
We observe that
$\ker\bigg(\Ext^1(\Z/p\Z,\Ga)\on{\fr-1}{\too}\Ext^1(\Z/p\Z,\Ga)\bigg)\simeq
\Z/p\Z$. It is possible to describe more explicitly
$\Ext^1(\Z/p\Z,\Z/p\Z)$. We recall that
$\Ext^1(\Z/p\Z,\Ga)=H^2_0(\Z/p\Z,\Ga)$ is freely generated as a
right $k$-module by the class of the  cocycle
$C_1=\frac{X^p+X^p-(X+Y)^p}{p}$. 

There is an isomorphism, induced by the maps of
\ref{prop:ext(Fp,Fp)},
$$
k/(\fr-1)(k) \times \Z/p\Z\too \Ext^1(\Z/p\Z,\Z/p\Z),
$$
 given by
$$
(a,b)\mapsto E_{a,b}.
$$
The extension $E_{a,b}$ is so defined: let $\bar{a}\in k$ a
lifting of $a$,
$$
E_{a,b}=\Sp(k[S_1,S_2]/(S_1^p-S_1,S_2^p-S_2-\bar{a} S_1))
$$

\begin{enumerate}
    \item law of multiplication
    \begin{align*}
    S_1\longmapsto &S_1\pt 1+1\pt S_1\\
    S_2\longmapsto &S_2\pt 1 + 1\pt S_2 + b \frac{S_1^p\pt 1+ 1\pt S_1^p-(S_1\pt 1+1\pt S_1)^p}{p}
    \end{align*}
    \item unit
    \begin{align*}
    &S_1\longmapsto 0\\
    &S_2\longmapsto 0
    \end{align*}
    \item inverse
    \begin{align*}
    &S_1\longmapsto -S_1\\
    &S_2\longmapsto -S_2
\end{align*}
\end{enumerate}
We remark that the extensions which are isomorphic to $\Z/p^2\Z$
as group schemes are the extensions $E_{0,b}$ with $b\neq 0$. By
\ref{ex:ext Z/pZ by Z/p } we have that any extension of $\Z/p\Z$
by $\Z/p\Z$ is given by $\clE^{(\lb_{(1)},\lb_{(1)};E_{p}(j \eta
S),j)}$.  We now study its reduction on the special fiber.

\begin{prop} For any $j\in \Z/p\Z$, $[\clE^{(\lb_{(1)},\lb_{(1)};E_{p}(j\eta S),j)}_k]=E_{0,j}\in
\Ext^1_k({\Z/p\Z},{\Z/p\Z})$. 

\end{prop}
\begin{proof}
As group schemes, $\clE^{(\lb_{(1)},\lb_{(1)};E_{p}(j \eta
S),j)}\simeq\Z/p^2\Z$, if $j\neq 0$, and
$\clE^{(\lb_{(1)},\lb_{(1)};1,0)}\simeq \Z/p\Z\times \Z/p\Z$
otherwise. In particular $\clE^{(\lb_{(1)},\lb_{(1)};E_{p}(j \eta
S),j)}_k$ has a  scheme-theoretic section. It is easy to see that
$\clE^{(\lb_{(1)},\lb_{(1)};E_{p}(j \eta S),j)}_k\simeq E_{0,b}$
with $$b=(-j\frac{\eta^p}{\lb_{(1)}(p-1)!} \mod \pi)=j,$$ since
$\frac{\eta^p}{\lb_{(1)}}\equiv
\frac{\lb_{(2)}^p}{\lb_{(1)}}\equiv 1\mod \pi$ and $(p-1)!\equiv
-1\mod \pi$ (Wilson Theorem). 
\end{proof}
\bibliographystyle{smfplain}
\bibliography{modelli21} \nocite{*}
\end{document}